\documentclass[11pt, a4paper]{amsart}

\usepackage[T1]{fontenc}
\usepackage{microtype}
\usepackage{parskip}

\usepackage{amsmath}
\usepackage{mathtools}
\usepackage{amssymb}
\usepackage{amsthm, bm}
\usepackage{amscd}
\usepackage{amsfonts}
\usepackage{MnSymbol}
\usepackage{enumerate}

\usepackage{stmaryrd}

\usepackage{extarrows}
\usepackage[bookmarks=true, colorlinks=true, linkcolor=blue!50!black,
citecolor=orange!50!black, urlcolor=orange!50!black, pdfencoding=unicode]{hyperref}
\usepackage{color}

\usepackage{tikz}									
\usetikzlibrary{matrix}
\usetikzlibrary{patterns}
\usetikzlibrary{matrix}
\usetikzlibrary{positioning}
\usetikzlibrary{decorations.pathmorphing}
\usetikzlibrary{cd}

\usepackage{fullpage}
\usepackage{cleveref}
\newtheorem*{theorem*}{Theorem}

\newtheorem{theorem}{Theorem}[section]
\newtheorem{lemma}[theorem]{Lemma}
\newtheorem{proposition}[theorem]{Proposition}
\newtheorem{corollary}[theorem]{Corollary} 

\theoremstyle{definition}

\newtheorem{example}[theorem]{Example}

\newtheorem{remark}[theorem]{Remark}

\newtheoremstyle{myitemstyle}						% flexible theorem style
	{}			%Space above
	{}			%Space below
	{}			%Body font
	{}			%indent amount
	{}			%Thm head font
	{.}			%Punkte nach thm head
	{ }			%Abstand nach thm head
	{}			%Thm head spec
\theoremstyle{myitemstyle}
\newtheorem{myitemthm}{}

\newcommand{\R}{\mathbb{R}}

\newcommand{\Z}{\mathbb{Z}}
\newcommand{\Q}{\mathbb{Q}}
\newcommand{\C}{\mathbb{C}}

\newcommand{\A}{\mathbb{A}}
\newcommand{\G}{\mathbb{G}}

\newcommand{\Res}{{\mathrm{Res}}}
\newcommand{\red}{{\mathrm{red}}}

\newcommand{\mbb}[1]{\mathbb{#1}}
\newcommand{\mc}[1]{{\mathcal{#1}}}

\DeclareMathOperator{\Hilb}{\mathrm{Hilb}}

\DeclareMathOperator{\Pic}{Pic}
\DeclareMathOperator{\Spec}{Spec}
\DeclareMathOperator{\Hom}{Hom}

\DeclareMathOperator{\codim}{codim}

\DeclareMathOperator{\id}{id}

\DeclareMathOperator{\Span}{Span}

\DeclareMathOperator{\GL}{GL}
\DeclareMathOperator{\Sp}{Sp}

\DeclareMathOperator{\SO}{SO}

\DeclareMathOperator{\Sym}{Sym}

\DeclareMathOperator{\SL}{SL}

\DeclareMathOperator{\ev}{ev}

\DeclareMathOperator{\ad}{ad}

\DeclareMathOperator{\conv}{conv}

\newcommand{\ab}{\mathrm{ab}}
\newcommand{\tf}{\mathrm{tf}}

\newcommand{\Dol}{{\mathrm{Dol}}}
\newcommand{\Betti}{{\mathrm{Betti}}}
\newcommand{\sd}{{\operatorname{sd}}}

\newcommand{\equ}{\ensuremath{\,=\,}}
\newcommand{\mbetti}{\mathcal{M}_B}

\title{P = W conjectures for varieties with nilpotent fundamental group}

\date{}

\author{Barbara Bolognese}
\address{Institut f\"ur Mathematik, Goethe--Universit\"at Frankfurt,
60325 Frankfurt am Main, Germany}
\email{bolognese@math.uni-frankfurt.de}

\author{Camilla Felisetti}
\address{Dipartimento di Scienze Fisiche, Informatiche e Matematiche, Università di Modena e Reggio Emilia.}
\email{camilla.felisetti@unimore.it}

\author{Andreas Gross}
\address{Institut f\"ur Mathematik, Goethe--Universit\"at Frankfurt,
60325 Frankfurt am Main, Germany}
\email{gross@math.uni-frankfurt.de}

\author{Alex K\"uronya}
\address{Institut f\"ur Mathematik, Goethe--Universit\"at Frankfurt,
60325 Frankfurt am Main, Germany}
\email{kuronya@math.uni-frankfurt.de}

\author{GyeongHyeon Nam}
\address{ Korea Institute for Advanced Study, 85 Hoegi-ro, Dongdaemun-gu, Seoul 02455, Republic of Korea}
\email{ghnam@kias.re.kr} 	

\author{Martin Ulirsch}
\address{Institut f\"ur Mathematik, Goethe--Universit\"at Frankfurt,
60325 Frankfurt am Main, Germany}
\email{ulirsch@math.uni-frankfurt.de}

\begin{document}

\begin{abstract}
 We prove the cohomological and the geometric $P=W$ conjectures for the moduli spaces of $G$-Higgs bundles and $G$-local systems on any smooth complex projective variety $X$ with nilpotent fundamental group, and any connected complex reductive group $G$.
To the best of our knowledge, our results provide the first instances of $P=W$ phenomena beyond curves and abelian varieties. 
Finally, we provide examples of our theorem for classical groups: as a byproduct of this analysis we investigate when symplectic resolutions of the corresponding character varieties  exist.
\end{abstract}

\maketitle

\setcounter{tocdepth}{1}
\tableofcontents

\section*{Introduction}
\renewcommand{\thetheorem}{\Alph{theorem}}

The structure of a complex algebraic variety decorates  the underlying topological space in many ways. A relatively new and intriguing aspect of this phenomenon is  the study of moduli spaces associated to the fundamental group of an algebraic variety. 

The topology of moduli spaces arising in non-abelian Hodge theory has been a central subject of study in algebraic geometry and representation theory since the foundational works of Hitchin, Donaldson, Corlette, and Simpson \cite{Cor88, Don87, Hit87, Simpson_nonabelianHodge, SimI, SimII}.

The theory is currently focused around extensions of the recently proven $P=W$ conjecture \cite{dCHM_Hilbert,dCHM_P=Wn=2, dCMS22, MS22, HMMS22, MSY25,FM22}, originally stated by de Cataldo--Hausel--Migliorini for certain moduli spaces over curves.  The main goal of our paper is to prove the cohomological and the geometric $P=W$ conjectures for moduli spaces of $G$-Higgs bundles and local systems over varieties with nilpotent fundamental group. 

Given a smooth complex projective variety $X$ and a complex reductive algebraic group $G$, non-abelian Hodge theory provides a canonical homeomorphism
\begin{equation}\label{eq:NAHT-isom}
\eta^G_X \colon M_{\Betti}^G(X) \xrightarrow{\quad \simeq \quad} M_{\Dol}^G(X)
\end{equation}
between the Betti moduli space $M_{\Betti}^G(X)$, the identity component of the character variety $\mathrm{Hom}(\pi_1(X), G) /\!\!/ G$,  and the Dolbeault moduli space $M_{\Dol}^G(X)$, which is the trivial component of the space of semistable $G$-Higgs bundles. While these two moduli spaces are far from isomorphic as algebraic varieties -- $M_{\Betti}^G(X)$ is affine, whereas when $X$ is a curve $M_{\Dol}^G(X)$ carries the structure of a quasi-projective non-compact integrable system -- the homeomorphism $\eta$ induces an isomorphism between their rational cohomology rings:
\[
H^*\big(M_{\Dol}^G(X), \mathbb{Q}\big) \cong H^*\big(M_{\Betti}^G(X), \mathbb{Q}\big) \ .
\]

When $X = C$ is a smooth projective curve of genus $g \ge 1$, the Dolbeault moduli space supports a proper morphism, the Hitchin fibration $h \colon M_{\Dol}^{G}(C) \to B$, which endows $M_{\Dol}^{G}(C)$ with the structure of an algebraically completely integrable system. By the decomposition theorem of Beilinson, Bernstein, Deligne, and Gabber \cite{BBD82}, the morphism $h$ equips the cohomology $H^*(M_{\Dol}^{G}(C), \mathbb{Q})$ with a geometric \emph{perverse Leray filtration} $P_\bullet H^*(M_{\Dol}^{G}(C), \mathbb{Q})$. On the Betti side, Deligne's mixed Hodge theory \cite{Del71, Del74} equips the cohomology of the affine variety $M_{\Betti}^{G}(C)$ with a \emph{weight filtration} $W_\bullet H^*(M_{\Betti}^{G}(C), \mathbb{Q})$.

Following earlier observations on \textquotedblleft curious Poincaré duality\textquotedblright\  by Hausel and Rodriguez-Villegas \cite{HauselRodriguezVillegas}, de Cataldo, Hausel, and Migliorini \cite{dCHM_P=Wn=2} formulated the celebrated \emph{$P=W$ conjecture}. 
It predicts that under the non-abelian Hodge isomorphism \eqref{eq:NAHT-isom}, the perverse Leray filtration $P_\bullet$ on the Dolbeault side corresponds to the weight filtration $W_\bullet$ on the Betti side, with a characteristic doubling of indices:
\begin{equation}\label{eq:P=W-curve}
P_k H^*\big(M_{\Dol}^{G}(C), \mathbb{Q}\big) = W_{2k} H^*\big(M_{\Betti}^{G}(C), \mathbb{Q}\big) = W_{2k+1} H^*\big(M_{\Betti}^{G}(C), \mathbb{Q}\big) \quad \text{for all } k \in \mathbb{Z}\  .
\end{equation}
The $P=W$ conjecture for curves has recently been established across various ranks and reductive groups through a series of major breakthroughs \cite{dCHM_P=Wn=2, dCMS_g=2pprime, dCMS22, MS22, HMMS22, MSY25,FM22} (see also \cite{Migliorini_survey,F25} for comprehensive surveys).

Parallel to the above, Katzarkov, Noll, Pandit, and Simpson \cite{KNPS15} introduced the \emph{geometric $P=W$ conjecture}, which is an effort to geometrically explain the surprising equality in the cohomological $P=W$ conjecture by comparing the geometry of neighborhoods at infinity of both moduli spaces that keep track of their quite different  compactifications. It is known by \cite[Theorem A]{MMS_geometricP=W} that, in the curve case, and under some technical assumptions, the geometric $P=W$ conjecture implies the top weight part of the cohomological $P=W$ conjecture: the relation between dual boundary complexes and top weight cohomology used there goes back to Payne \cite{Payne13}. 
The geometric  $P=W$ conjecture asserts that the dual boundary complex $\Delta(\mbetti)$ associated with a divisorial log terminal (dlt) log Calabi--Yau compactification of $M_{\Betti}^G(X)$ is homotopy equivalent (or PL-homeomorphic) to the boundary sphere quotient of the Hitchin base, aligning the asymptotic geometry at infinity of both moduli spaces \cite{MMS_geometricP=W}. For instance, Simpson showed that the dual boundary complex of the $\SL_2$-character variety of a punctured sphere is homotopy equivalent to a sphere of the expected dimension \cite{Simpson16}, and the conjecture was proved for curves of genus one and $G=\GL_n$ or $\SL_n$ in \cite{MMS_geometricP=W}.

\subsection*{$P=W$ phenomena beyond curves}

It is natural to ask what happens when we replace  curves in the above discussion  by a higher dimensional smooth projective variety $X$. First, new difficulties arise already in the definition of the objects involved. In fact, when $\dim X\geq 2$, a Higgs field $\theta$ must satisfy a certain integrability condition (which is vacuous for curves), and, as a consequence, the Hitchin map 
\[
h^G_X\colon M_{\Dol}^G(X)\rightarrow \mathcal{A}_X^G\coloneqq\bigoplus_i H^0(X, S^{e_i}\Omega^1_X)
\]
is no longer surjective. Chen and Ng$\hat{\mathrm{o}}$ first addressed this problem in \cite{CN20}, and went on to  show that $h_G^X$ factors through a \textit{spectral data morphism} $\sd^G_X\colon M_{\Dol}^G(X)\rightarrow B_X^G$, whose target takes the integrability condition into account. They conjectured that the spectral data morphism is surjective; this has been verified in several cases, see for instance \cite{PatelWeissmann26, SheshmaniWangXia26}, while counterexamples have recently been constructed in \cite{HeLiuZhang26}.

Prior to this work, higher dimensional $P=W$ phenomena in non-abelian Hodge theory were only known for abelian varieties and $G=\GL_n$ (for $P=W$ phenomena of a different nature, concerning Lagrangian fibrations of compact hyper-K\"ahler manifolds, see \cite{ShenYin22, HLSY21,FSY22}). In \cite{BKU23}, the authors established the cohomological and geometric $P=W$ conjectures. It is worth underlining that curves of positive genus and abelian varieties are both aspherical, so their homotopy type is entirely determined by their fundamental group.

In this paper we prove the cohomological and the geometric $P=W$ conjectures for an arbitrary connected reductive group $G$, and for every smooth projective variety $X$ with nilpotent fundamental group. 
This class of varieties is very large: it contains, for example, varieties of general type in every dimension $n\geq 2$ which are not aspherical (cf. Corollary B). As far as we know, this is the first case in which $P=W$ conjectures are proven outside the context of curves and abelian varieties. 

While the main $P=W$ results are established for varieties with nilpotent fundamental group, several of our constructions and results hold for an arbitrary smooth projective variety $X$, laying the foundations for a systematic study of $P=W$ phenomena in higher dimension. 

We use the nilpotency of $\pi_1(X)$ only once, in Lemma \ref{lem:unipotent on identity component}: together with \Cref{lem:unipotent on identity component}, it implies that every semisimple $G$-representation of $\pi_1(X)$ in the connected component of the trivial one factors through the free abelian group $H_1(X,\mathbb{Z})/\mathrm{tors}$, and can therefore be conjugated into a maximal torus~$T$.
Up to normalizations, this allows for describing our moduli spaces as finite quotients of the corresponding moduli spaces where $G$ is replaced by $T$. In this way, our argument makes visible and explicit the precise role played by the nilpotent nature of the fundamental group. 

\subsection*{Spectral data morphisms and cohomological $P=W$}
Fix a maximal torus $T\subseteq G$ and let $N$ be its cocharacter lattice, and $W$ the associated Weyl group. \\
For the torus $T$, one observes  that the spectral data morphism is  the projection
\[
\sd^T_X\colon M_{\Dol}^T(X)\cong \left(\Pic^0(X)\times H^0(X,\Omega^1_X)\right)\otimes_{\Z}N\rightarrow B^T_X= H^0(X,\Omega^1_X)\otimes_{\Z}N\ .
\]
If $\pi_1(X)$ is nilpotent, we show that $\widetilde M_{\Dol}^G(X)\cong M_{\Dol}^T(X)/W$ is the normalization of $M_{\Dol}^G(X)$, and that it is homeomorphic to $M_\Dol^G(X)$. The map $\sd^T_X$ descends to a proper surjective morphism $\widetilde{\sd}^G_X$ from this normalization to $\widetilde{B}^G_X\coloneqq B^T_X/W$. 
We establish that $\sd^G_X,\widetilde{\sd}_X^G,$ and $h^G_X$ all induce the same perverse filtration on the cohomology of $M_{\Dol}^G(X)$, which we can therefore  denote by $P_{\bullet}H^*(M_{\Dol}^G(X))$ without any ambiguity.\footnote{It is worth noticing that, while the result is established for $\pi_1(X)$ nilpotent, the two main ingredients of the proof -- which are a functoriality lemma on the perverse filtration and the finiteness of the composition between the map $\beta^G_X\colon \widetilde{B}^G_X\to B^G_X$ induced by the inclusion $T\subseteq G$ and the natural morphism $\iota^G_X\colon B^G_X\rightarrow \mathcal{A}^G_X$ -- are proved for arbitrary $X$.}

Armed with this, we are now in a position to state our first main result,  the cohomological $P=W$ conjecture.

\begin{theorem}
\label{thm:P=W}
Let $X$ be a smooth projective variety with nilpotent group and let $G$ be a complex connected reductive group.
Then the  non-abelian Hodge correspondence 
\[
    \eta_X^G\colon M^G_\Betti(X)\to M^G_\Dol(X)
\]
induces an isomorphism
\[
    P_k H^*(M^G_\Dol(X),\Q) \xrightarrow{\cong}
    W_{2k} H^*(M^G_\Betti(X),\Q) \,=\,
    W_{2k+1} H^*(M^G_\Betti(X),\Q)\ .
\]
%where the perverse filtration is taken with respect to the spectral data morphism.
\end{theorem}

We note that, for example  by the Lefschetz hyperplane theorem, varieties with an abelian fundamental group are abundant, and we obtain the following corollary.

\begin{corollary}
	\label{cor:intro general type}
	For every $n\geq2$ there exist smooth complex projective varieties $X$ of general type and of dimension $n$, which are not aspherical, such that the conclusion of Theorem \ref{thm:P=W} holds for every connected complex reductive group $G$. For instance, one can take for $X$ a general complete intersection of very ample divisors in an abelian variety of dimension larger than $n$.
\end{corollary}

Next, we  turn to the geometric $P=W$ conjecture. In this setting the role of the Hitchin base is played by $\widetilde{B}^G_X$, and the role of the sphere at infinity by $\widetilde{B}^G_X\setminus\{0\}$. The latter  is homotopy equivalent to $S^{2gr-1}/W$, where $2g$ is the rank of the abelianization of $\pi_1(X)$ and $r$ is the rank of $G$. 

\begin{theorem}
\label{thm:geometric P=W}
Let $X$ be a smooth projective variety with nilpotent fundamental group and let $G$ be a complex connected reductive group. 
There exists a dlt model  for the normalization $\widetilde M^G_\Betti(X)$ of $M^G_\Betti(X)$ admitting a neighborhood  $N^G_\Betti(X)$ at infinity, a  dual complex $\Delta^G_X$, and an evaluation map $\ev\colon N^G_\Betti(X)\to \Delta^G_X$, such that there exists a homotopy equivalence $\chi\colon \Delta_X^G \to \widetilde B^G_X\setminus \{0\}$ making the diagram
\[
\begin{tikzcd}
    N^G_\Betti(X) 
    \arrow[r,"\eta_X^G"]
    \arrow[d,"\ev"]&
        N^G_\Dol(X)
        \arrow[d,"\widetilde\sd^G"]\\
    \Delta^G_X
    \arrow[r,"\chi"]&
        \widetilde B^G_X\setminus\{0\}
\end{tikzcd}
\]
commute up to homotopy. Here, $N^G_\Dol(X)=\eta_X^G(N^G_\Betti(X))$.
\end{theorem}

We note that if $X$ is an elliptic curve, then $\widetilde B^G_X\setminus \{0\}$ is homotopy equivalent to a sphere, making the theorem consistent with the geometric $P=W$ conjecture as stated in \cite{KNPS15}. In particular, for $G=\GL_n$ or $\SL_n$, Theorem \ref{thm:geometric P=W} recovers the genus one case of the geometric $P=W$ conjecture proved in \cite{MMS_geometricP=W}, since in these cases $\widetilde M^G_\Betti(X)=M^G_\Betti(X)$ and $\widetilde B^G_X$ is the Hitchin base. For $g\geq 2$, the base $\widetilde B^G_X\setminus \{0\}$ is not homotopy equivalent to a sphere in general \cite{BKU23}.
Also, for classical groups, the moduli space $M_{\Betti}^G(X)$ is known to be normal \cite[Theorem~2.1(3)]{Sikora}, so that we acquire an isomorphism $M_{\Betti}^G(X)\cong \widetilde M_{\Betti}^G(X)$.
 
We end the paper by making the cohomological results explicit in the case of classical groups of type $B_n,C_n$, and $D_n$. We show that the moduli spaces for $\Sp_{2n}$ and $\SO_{2n+1}$ are symmetric products of a singular variety, while  those  for $\SO_{2n}$ are double covers of the former. 
In our concrete cases, the study of the spectral data morphism boils down to understanding the map from the singular variety to an affine space, and  considering the induced map on the symmetric products. As a byproduct of this analysis, we determine when these moduli spaces admit a symplectic resolution. 

\begin{corollary}
	\label{cor:intro symplectic}
	Let $X$ be a smooth complex projective variety with a nilpotent fundamental group whose abelianization has rank $2g\geq2$, and let $G$ be a connected complex reductive group which is not a torus. Then the normalizations of $M^G_\Dol(X)$ and $M^G_\Betti(X)$ have symplectic singularities; in particular, they have rational Gorenstein singularities. Moreover,
	\begin{enumerate}
		\item if $g\geq2$, neither of them admits a symplectic resolution;
		\item if $g=1$ and $G$ is one of $\GL_n$, $\SL_n$, $\Sp_{2n}$ and $\SO_{2n+1}$, both of them admit a symplectic resolution, given by a Hilbert scheme of points on a surface or by a fiber of it;
		\item if $g=1$ and the root system of $G$ has an irreducible component which is not of type $A$, $B$ or $C$, neither of them admits a projective symplectic resolution.
	\end{enumerate}
\end{corollary}

For $g=1$ the existence of symplectic resolutions has been studied for groups of type $A$ by Bellamy and Schedler \cite{BellamySchedler23}, and for all almost simple groups by Shu \cite{Shu24}, who showed that the answer depends on the isogeny type of $G$. For character varieties of free abelian groups, results closely related to Corollary \ref{cor:intro symplectic} were obtained, while this paper was being written, in \cite{FGPZ26}. We compare these results with ours in Section \ref{sec:symplectic-resolutions}.

\subsection*{Structure of the paper}
In Section \ref{sec:functoriality}, we study perverse filtrations under finite group actions and finite morphisms, proving the functoriality properties of the perverse filtration that we need later. In Section \ref{sec:NAH}, we recall the Betti and Dolbeault moduli spaces and the non-abelian Hodge correspondence for varieties with nilpotent fundamental group. We then recall the construction of the spectral data morphism for arbitrary $X$, compare the perverse filtrations associated with the spectral data and Hitchin morphisms, and prove Theorem \ref{thm:P=W} and Corollary \ref{cor:intro general type}. Section \ref{sec:geometric P=W} is devoted to the geometric $P=W$ conjecture: after recalling evaluation maps and dual complexes, we treat projective toric varieties and their finite quotients, and we prove Theorem \ref{thm:geometric P=W}. In Section \ref{sec:pwclassicalgroups}, we make our results explicit for the classical groups, we study the singularities of the moduli spaces and prove Corollary \ref{cor:intro symplectic}.

\subsection*{Acknowledgements}
Helpful and motivating discussions with Jakob Stix and Torsten Wedhorn were much appreciated. 
B.~B., A.~G., and A.~K. have  received funding from the Deutsche For\-schungs\-gemeinschaft (DFG, German Research Foundation) TRR 326 \emph{Geometry and Arithmetic of Uniformized Structures}, project number 444845124. In addition, A.~G. has been supported by the Marie-Sk\l{}odowska-Curie-Stipendium Hessen (as part of the HESSEN HORIZON initiative), and A~K. by  the ANR-DFG project “Positivity
on K-trivial varieties” (project number ANR-23-CE40-0026).
C.\ F.\ was supported by IN$\delta$AM GNSAGA and by the project “Discrete Methods in Combinatorial Geometry and Geometric Topology” of the University of Modena and Reggio Emilia.
G.\ N.\ was supported by Oscar Kivinen's Väisälä project grant of the Finnish Academy of Science and Letters and a KIAS Individual Grant (HP111501) at Korea Institute for Advanced Study.

\renewcommand{\thetheorem}{\arabic{section}.\arabic{theorem}}

%%%%%%%%%%%%%%%%%%%%%%%%%%%%%%%%%%%%%%%%%%%%%%%%%%%
%% SECTION 1
%%%%%%%%%%%%%%%%%%%%%%%%%%%%%%%%%%%%%%%%%%%%%%%%%%%
\section{Functoriality of Perverse Filtrations under Finite Group Actions}\label{sec:functoriality}

We collect the necessary technical results about perverse filtrations. 
Let $X$ and $Y$ be complex quasi-projective varieties, and let $D^b_c(Y)$ denote the bounded derived category of constructible sheaves on $Y$ with rational coefficients. All hypercohomology groups are taken with coefficients in $\mathbb{Q}$. In what follows, we will recall only the definitions and concepts which are needed in this paper; for the general theory of perverse sheaves we refer the reader to \cite{BBD82, dCM09}.
	
Given a proper algebraic map $f\colon X \to Y$, one defines the \textit{perverse Leray filtration} on $H^*(X, \mathbb{Q})$ via
\begin{equation}\label{def:perversefiltration}
	P_p^f H^d(X, \mathbb{Q}) \coloneqq \mathrm{Im} \left( H^{d - \dim X + r(f)}\left( Y, {}^{\mathfrak{p}}\tau_{\le p} R f_* \mathbb{Q}_X [\dim X - r(f)] \right) \longrightarrow H^d(X, \mathbb{Q}) \right)\ ,
\end{equation}
where $r(f) \coloneqq \dim(X \times_Y X) - \dim X$ is called the \textit{defect of semismallness} of $f$ and ${}^{\mathfrak{p}}\tau_{\le p}$ denotes the \textit{perverse truncation functor} on $D^b_c(Y)$. Usually, the map $f$ is clear from context, in which case we drop it from the notation.
	
\subsection{Perverse filtration and finite group actions}
	
Let $W$ be a finite group acting equivariantly on a proper morphism $f \colon X \to Y$, inducing a commutative diagram of quotient varieties
\[
\begin{tikzcd}
    X \arrow[r, "q_X"] \arrow[d, "f"'] & X/W \arrow[d, "f_W"] \\
    Y \arrow[r, "q_Y"']               & Y/W\ .
\end{tikzcd}
\]
The following result is probably well known to experts, nevertheless,  we provide a short proof for the lack of a suitable reference.  
	
\begin{proposition}\label{prop:quotient-filtration}
		Let $f \colon X \to Y$ be a $W$-equivariant proper morphism between quasi-projective varieties, and let $f_W \colon X/W \to Y/W$ be the induced morphism on quotients. Under the canonical identification $H^*(X/W, \mathbb{Q}) \cong \left( H^*(X, \mathbb{Q}) \right)^W$, the perverse filtration on $H^*(X/W, \mathbb{Q})$ associated with $f_W$ satisfies
		\[
		P_p H^*(X/W, \mathbb{Q}) = \left( P_p H^*(X, \mathbb{Q}) \right)^W.
		\]
\end{proposition}
\begin{proof}
By definition, the perverse filtration on $H^*(X/W, \mathbb{Q})$  with respect to $f_W$ is given by the image of the canonical hypercohomology morphism induced by the perverse truncation functor
\[
H^*\left( Y/W, {}^{\mathfrak{p}}\tau_{\le p} Rf_{W*} \mathbb{Q}_{X/W}[\dim X/W-r(f_W)] \right) \longrightarrow H^*(X/W, \mathbb{Q})\ .
\]
Since the quotient map $q_Y \colon Y \to Y/W$ is finite, the derived pushforward $q_{Y*} = Rq_{Y*}$ is exact with respect to the perverse $t$-structure \cite[Corollaire 2.2.6]{BBD82}. The $W$-action on $X$ induces automorphisms of $R f_* \mathbb{Q}_X$ in $D^b_c(Y)$. Via the identity $q_Y \circ w = q_Y$, these automorphisms descend to endomorphisms of $q_{Y*} R f_* \mathbb{Q}_X$ in $D^b_c(Y/W)$, and allow us to define the projector
\[
p_W \coloneqq \frac{1}{|W|} \sum_{w \in W} q_{Y*}(w^\#) \;\in\; \operatorname{End}_{D^b_c(Y/W)}\big( q_{Y*} R f_* \mathbb{Q}_X \big)\ .
\]
Working over $\mathbb{Q}$, the operator $p_W$ is an idempotent in the idempotent-complete category $D^b_c(Y/W)$, yielding the canonical direct sum decomposition
\[
q_{Y*} R f_* \mathbb{Q}_X \cong Rf_{W*} \mathbb{Q}_{X/W} \oplus \mathcal{K}^\perp, \quad \text{with} \quad Rf_{W*} \mathbb{Q}_{X/W} = \operatorname{Im}(p_W)\ .
\]
Since ${}^{\mathfrak{p}}\tau_{\le p} \colon D^b_c(Y/W) \to D^b_c(Y/W)$ is an additive functor commuting with the $W$-action, it commutes with the projector $p_W$. Consequently, perverse truncation commutes with taking invariants:
\[
{}^{\mathfrak{p}}\tau_{\le p} Rf_{W*} \mathbb{Q}_{X/W} \cong \left( {}^{\mathfrak{p}}\tau_{\le p} q_{Y*} R f_* \mathbb{Q}_X \right)^W.
\]
In particular, taking the image of the evaluation map commutes with taking invariants, which yields
\begin{align*}
P_pH^*(X/W,\mathbb{Q})&=\mathrm{Im}\left( H^{*-\dim X+r(f_W)}\!\left(Y/W, {}^{\mathfrak{p}}\tau_{\le p} Rf_{W*} \mathbb{Q}_{X/W}[\dim X-r(f_W)] \right) \to H^*\!(X/W, \mathbb{Q}) \!\right) \\
        &= \mathrm{Im}\left( H^{*-\dim X+r(f)}\left( Y/W, \left( {}^{\mathfrak{p}}\tau_{\le p} q_{Y*}R f_* \mathbb{Q}_X[\dim X-r(f)] \right)^W \right) \rightarrow \left( H^*(X, \mathbb{Q}) \right)^W \right) \\
            &= \left( \mathrm{Im}\left( H^{*-\dim X+r(f)}\left( Y, {}^{\mathfrak{p}}\tau_{\le p} R f_* \mathbb{Q}_X [\dim X-r(f)]\right) \rightarrow H^*(X, \mathbb{Q}) \right) \right)^W \\
                &= \left( P_p H^*(X, \mathbb{Q}) \right)^W \ .
\end{align*}
This completes the proof.
\end{proof}
   
More generally, we can state the following technical result, which is of independent interest, and describes the behavior of the perverse filtration of a map under the composition with finite morphisms on the left and on the right.

\begin{lemma}\label{lem:perversecompostion}
Let $f\colon X\rightarrow Y$ be a proper algebraic map of complex quasi-projective varieties. 
\begin{enumerate}[(i)]
\item If $g\colon Y\rightarrow Z$ is a finite morphism to an algebraic variety $Z$ then 
\[ 
    P_p^{g\circ f}H^*(X,\mathbb{Q})=P_p^fH^*(X,\mathbb{Q}) \text{ for all }p\in \mathbb{Z}\ .
\]
\item If $\nu\colon\widetilde{X}\rightarrow X$ is a finite morphism from an algebraic variety $\widetilde X$ which is bijective on points, then the pullback $\nu^*\colon H^*(X,\mathbb{Q})\rightarrow H^*(\widetilde{X},\Q)$ is an isomorphism that is compatible with the perverse filtration, i.e.
\[ 
    P_p^{f\circ \nu}H^*(\widetilde{X},\mathbb{Q})=\nu^*P_p^fH^*(X,\Q)\text{ for all }p\in \mathbb{Z}\ .
\]
\end{enumerate}
\end{lemma}
\begin{proof}
We start by proving (i). The first thing we must show is that the respective perverse filtrations of $f$ and $g\circ f$ have the same range, that is, $r(g\circ f)=r(f)$. On the one hand, since $g$ is separated, the fiber product $X\times_YX$ over $Y$ is a closed subscheme of the fiber product $X\times_ZX$ over $Z$, hence $r(f)\leq r(g\circ f)$. Conversely, we set $Y_m\subseteq Y$ to be the set of $y\in Y$ satisfying $\dim f^{-1}(y)=m$. It is easy to see that $\dim X\times_Y X\geq \dim Y_m+2m$.
In addition, since $g$ is finite,  we have that both projections
\[ 
    Y_m\leftarrow Y_m\times_ZY_n\rightarrow Y_n
\]
are quasi-finite for any $m,n$, hence $\dim Y_m\times_ZY_n\leq \min\{\dim Y_m, \dim Y_n\}$.
If we stratify $X\times_ZX$ with respect to the dimensions of the fibers of $f$ we get
\begin{align*}
\dim(X\times_ZX) &\leq \max_{m,n}\left\lbrace\min\{\dim Y_m, \dim Y_n\}+m+n\right\rbrace\\
 &\leq \frac{1}{2}\max_{m,n}\left\lbrace \dim Y_m +2m +\dim Y_n+2n\right\rbrace\leq \dim (X\times_YX)\ .
\end{align*}
Subtracting $\dim X$ from both sides we get  $r(g\circ f)\leq r(f)$.
This implies, in particular, that the defects of semismallness of $f$ and $g\circ f$ coincide, and their perverse filtrations are defined with the same shift in \eqref{def:perversefiltration}.

As before, since $g$ is finite, we obtain that the functor $g_*$ is perverse t-exact, therefore it commutes with the perverse truncation functors $^{\mathfrak p}\tau_{\leq p}$. Thus, we have 
\begin{align*}
P^{g\circ f}_pH^*(X,\mathbb{Q})&=\mathrm{Im}\left( H^{*-\dim X+r(g\circ f)}\left( Z, {}^{\mathfrak{p}}\tau_{\leq p} g_*R f_* \mathbb{Q}_X[\dim X-r(g\circ f)] \right) \longrightarrow H^*(X, \mathbb{Q}) \right)\\
&=\mathrm{Im}\left( H^{*-\dim X+r(g\circ f)}\left( Z, g_* {}^{\mathfrak{p}}\tau_{\le p} R f_* \mathbb{Q}_X [\dim X-r(g\circ f)]\right) \longrightarrow H^*(X, \mathbb{Q}) \right)\\
&=\mathrm{Im}\left( H^{*-\dim X+r(f)}\left(Y, {}^{\mathfrak{p}}\tau_{\le p} R f_* \mathbb{Q}_X [\dim X-r(f)]\right) \longrightarrow H^*(X, \mathbb{Q}) \right)\\
&=P_p^fH^*(X,\mathbb{Q})\ .
\end{align*}
This concludes the proof of (i). 
For (ii), it is straightforward to see that, since $\nu$ is finite and bijective, the morphism $\mathbb{Q}_X\rightarrow \nu_*\Q_{\widetilde{X}}$ is an isomorphism. This implies that $\nu_*$ induces an isomorphism also at the level of complexes
\[
    Rf_*\mathbb{Q}_X\cong Rf_*\nu_*\mathbb{Q}_{\widetilde{X}}=R(f\circ \nu)_*\mathbb{Q}_{\widetilde{X}}\ .
\]
The above isomorphism implies that $\nu^*$ preserves the perverse filtration in cohomology, provided that we check that the shifts in the definition \eqref{def:perversefiltration} coincide. However, this is true since the hypotheses on $\nu$ imply that $r(f\circ\nu)= \dim (\widetilde X\times_Y\widetilde X)-\dim \widetilde X=\dim (X\times_YX)-\dim X=r(f)$.
\end{proof}

\begin{remark}\label{rem:pervdimconv}
Lemma \ref{lem:perversecompostion} relies on the fact that the normalization of the perverse filtration that we adopt only depends on the morphism itself and not on the dimension of the target. This fact will be crucial in Section 2, where we are using this lemma to prove that, when $\pi_1(X)$ is nilpotent, the perverse filtration of the Hitchin map coincides with that of the various spectral data morphisms (cf.\  Proposition \ref{prop:allperverseequal}).
\end{remark}

%%%%%%%%%%%%%%%%%%%%%%%%%%%%%%%%%%%%%%%%%%%%%%%%%%%
%% SECTION 2
%%%%%%%%%%%%%%%%%%%%%%%%%%%%%%%%%%%%%%%%%%%%%%%%%%%
\section{Non-abelian Hodge correspondence and cohomological $P=W$}\label{sec:NAH}

This section is devoted to the proof of the cohomological  $P=W$ conjecture for moduli spaces of $G$-Higgs bundles bundles over varieties with nilpotent fundamental group. 
We will see that the nilpotency assumption can be weakened: it suffices that $\pi_1(X)$ becomes nilpotent after dividing out the subgroup generated by its elements of finite order (cf. Proposition \ref{prop:reduction to T Betti}), which yields examples with non-nilpotent fundamental group, see Example \ref{ex:non-nilpotent}.

First, we recall the non-abelian Hodge correspondence on smooth projective varieties \cite{Simpson_nonabelianHodge} along with its protagonists.  Throughout this section, let $G$ be a connected reductive group with maximal torus $T$. Let $N$ denote the cocharacter lattice of $T$, and let $W$ denote the Weyl group of $G$. Moreover, let $X$ be a smooth projective  variety.

\subsection{The Betti moduli space}

Let $\Gamma$ be a finitely presented group. Then the set of representations $\Hom(\Gamma,G)$ carries the structure of an algebraic variety and the GIT quotient 
\[
\mc X_G(\Gamma)\coloneqq \Hom(\Gamma,G)/\!\!/ G\ ,
\]
where $G$ acts by conjugation, is the \emph{character variety}. We denote by $\mc X_G(\Gamma)_0$ the connected component of the trivial representation. 

Now let $X$ be a smooth projective variety over the complex numbers.
We define the \emph{Betti moduli space} $M^G_\Betti(X)$ as $\mc X_G(\pi_1(X))_{0,\red}$, that is $\mc X_G(\pi_1(X))_0$ equipped with its induced reduced scheme structure.

Any morphism $H\to G$ of reductive groups induces a morphism $\mc X_H(\Gamma)\to \mc X_G(\Gamma)$. In particular, 
the inclusion $T\to G$ induces a morphism
\(
    M_\Betti^T(X)\to  M^G_\Betti(X)
\). This  is $W$-invariant and so it yields a  morphism
\[
 \Phi\colon    M^{T}_\Betti(X)/W \to M^G_\Betti(X)\ .
\]

The space $M_\Betti^T(X)$, or more generally $\mc X_T(\Gamma)$, is readily understood. Because $T$ is abelian, we have
\[
    \mc X_T(\Gamma)=\Hom(\Gamma^\ab,T)=\Hom(\Gamma^\ab,\G_m)\otimes_\Z N\ ,
\]
where $\Gamma^\ab$ denotes the abelianization of $\Gamma$. If $\Gamma^\ab$ has torsion, the space $\mc X_T(\Gamma)$ is not connected. For the component of the trivial representation, we have
\[
    \mc X_T(\Gamma)_0 = \Hom(\Gamma^\ab_\tf,\G_m)\otimes_\Z N \ ,
\]
where $\Gamma^\ab_\tf=\Gamma^\ab/\mathrm{Tor}(\Gamma^\ab)$ is the torsion free part of $\Gamma^\ab$.
In particular, $\mc X_T(\Gamma)_0$ is an algebraic torus.

Our goal in this section is to prove the following result:
\begin{proposition}
\label{prop:reduction to T Betti}
    Suppose that the quotient of $\pi_1(X)$ by its normal subgroup generated by all elements of finite order is nilpotent. Then the natural morphism 
    \[
    \Phi\colon M^{T}_\Betti(X)/W \to M^G_\Betti(X)
    \] 
    is bijective and finite, and in particular a homeomorphism. Moreover, if $\widetilde M^G_\Betti(X)$ denotes the normalization of $M^G_\Betti(X)$, then $\Phi$ factors as an isomorphism $M^T_\Betti(X)/W\to \widetilde M^G_\Betti(X)$ followed by a homeomorphism $\widetilde M^G_\Betti(X)\to M^G_\Betti(X)$.
\end{proposition}

Before we can prove the proposition, we need several auxiliary results about character varieties. 

\begin{lemma}
\label{lem:unipotent on identity component}
Let $\Gamma$ be a finitely presented group and let $s\in \Gamma$. Further suppose that $\psi(s)$ has finite order for every finite-dimensional irreducible representation $\psi\colon \Gamma\to \GL(V)$. 
Then $\rho(s)$ is unipotent for every $\rho\in \Hom(\Gamma,G)_0$. If, moreover, $\rho$ is semisimple, that is, if its $G$-orbit is closed, then $\rho(s)=1$. 
\end{lemma}

\begin{proof}
Fix a faithful representation $\iota\colon G\hookrightarrow\GL(V)$, and let $n=\dim V$. Let $\psi\in \Hom(\Gamma,G)$. The eigenvalues of $\iota\circ \psi(s)$ are its eigenvalues on the composition factors of the representation  $\iota\circ \psi$ of $\Gamma$. These composition factors are irreducible representations of $\Gamma$, on which $s$ acts with finite order by assumption, and therefore all the eigenvalues of $\iota\circ\psi(s)$ are roots of unity.

Now, consider the map
\[
    \varphi\colon\Hom(\Gamma,G)_0\longrightarrow\C^{n}\ ,
\]
which sends $\psi$ to the coefficients of the characteristic polynomial of $\iota\circ\psi(s)$. It takes values in the countable set of monic polynomials all of whose roots are roots of unity. On the other hand, its image is constructible by Chevalley's theorem, and a countable constructible subset of $\C^n$ is finite. 
As a result, $\Hom(\Gamma,G)_0$ is the disjoint union of finitely many fibers of $\varphi$, which are closed, and hence also open. Since $\Hom(\Gamma,G)_0$ is connected, $\varphi$ is constant.
Evaluating at the trivial representation, we find that the characteristic polynomial of $\iota\circ\rho(s)$ is $(x-1)^{n}$, that is that $\rho(s)$ is unipotent for every $\rho\in \Hom(\Gamma,G)_0$.

Now assume that $\rho\in \Hom(\Gamma,G)_0$ is semisimple. Then the Zariski closure of $\rho(\Gamma)$ is reductive by \cite[Theorem 3.6]{Richardson}, so that $\iota\circ \rho$ is a direct sum of irreducible representations of $\Gamma$. By assumption, $s$ acts with finite order on each of them, thus $\iota\circ\rho(s)$ is unipotent and of finite order, hence trivial. Since $\iota$ is faithful, it follows that $\rho(s)=1$.
\end{proof}

\begin{lemma}
\label{lem:quotient by finite-order-in-rep-elements}
Let $\Gamma$ be a finitely presented group and $S\subseteq \Gamma$ a subset such that $\psi(s)$ has finite order for every $s\in S$ and every finite-dimensional irreducible representation $\psi\colon \Gamma\to \GL(V)$ of $\Gamma$. Further, assume that $\Gamma/\llangle S\rrangle$ is finitely presented, where $\llangle S\rrangle$ denotes the normal subgroup generated by $S$. Then the morphism
\[
    j\colon \mc X_G(\Gamma/\llangle S\rrangle)_{0,\red}\to \mc X_G(\Gamma)_{0,\red}
\]
induced by the quotient map $p\colon\Gamma\to \Gamma/\llangle S\rrangle$ is an isomorphism.
\end{lemma}
\begin{proof}
The surjection $p$ identifies $\Hom(\Gamma/\llangle S\rrangle,G)$ with the closed $G$-stable subscheme of $\Hom(\Gamma,G)$ of representations which are trivial on $S$. Since $G$ is reductive, the induced map on rings of invariants is surjective, so that the induced morphism $\mc X_G(\Gamma/\llangle S\rrangle)\to\mc X_G(\Gamma)$, which we still denote by $j$, is a closed immersion.

As $G$ is connected, the connected components of $\Hom(\Gamma,G)$ are $G$-stable, and therefore $\mc X_G(\Gamma)_0=\Hom(\Gamma,G)_0/\!/G$. In particular, every point of $\mc X_G(\Gamma)_0$ is represented by a semisimple representation $\rho\in\Hom(\Gamma,G)_0$.

By Lemma \ref{lem:unipotent on identity component}, such a $\rho$ factors through $p$, so that $\mc X_G(\Gamma)_0$ is contained in the image of $j$. Since $j$ is a homeomorphism onto its closed image and $\mc X_G(\Gamma)_0$ is connected and contains the class of the trivial representation, it follows that $\mc X_G(\Gamma)_0\subseteq j(\mc X_G(\Gamma/\llangle S\rrangle)_0)$. Conversely, $j(\mc X_G(\Gamma/\llangle S\rrangle)_0)$ is connected and contains the class of the trivial representation, hence is contained in $\mc X_G(\Gamma)_0$. As a result, $j$ restricts to a surjective closed immersion of $\mc X_G(\Gamma/\llangle S\rrangle)_0$ into $\mc X_G(\Gamma)_0$, and hence to an isomorphism of the induced reduced schemes.
\end{proof}

\begin{lemma}
\label{lem:roots of unity}
Let $\Gamma$ be nilpotent, let $\psi\colon \Gamma\to \GL(V)$ be an irreducible representation on a finite-dimensional complex vector space $V$, and let $s\in [\Gamma,\Gamma]$. Then $\psi(s)$ has finite order.
\end{lemma}
\begin{proof}

Let $H=\overline{\psi(\Gamma)}\subseteq\GL(V)$ be the Zariski closure of $\psi(\Gamma)$. By \cite[\S 4.3]{Waterhouse79}, $H$ is again nilpotent, since $[\overline A,\overline B]\subseteq\overline{[A,B]}$ for all subgroups $A,B\subseteq\GL(V)$. In particular, its identity component $H^{0}$ is connected and solvable, thus by the Lie--Kolchin theorem \cite[\S 10.2]{Waterhouse79} we can choose a basis of $V$ in which all elements of $H^0$ are upper triangular. 
The unipotent elements of $H^0$ then form a subgroup $H^0_u$, which is given by the kernel of the homomorphism which sends an element of $H^0$ to its diagonal. Since $H^0$ is normal in $H$ and conjugation preserves unipotent elements, $H^0_u$ is normal in $H$ and the subspace $V^{H^0_u}\subseteq V$ of $H^0_u$-invariants is stable under the action of $H$.
By construction $V^{H^0_u}$ contains the first vector of the chosen basis, hence, since $V$ is irreducible, we have that $V^{H^0_u}=V$ and $H^0_u=\{1\}$. This implies that $H^0$ embeds in the diagonal torus, and thus that it is commutative and it consists of semisimple elements. Since $H^0$ is connected, it is a torus.

We now want to show that $H^0$ is central in $H$. Conjugation defines a homomorphism 
\[
    H\to\mathrm{Aut}(H^0)=\GL(X_*(H^0)),\; h\mapsto\alpha_h,
\]
which is trivial on $H^0$ and hence has finite image.
If $H$ has nilpotency class at most $c$, then the $c$-fold commutator $[h,[h,\dots,[h,t]\dots]]$ is trivial for all $h\in H$ and $t\in H^0$. In additive notation on $H^0=X_*(H^0)\otimes_\Z\C^*$ this commutator equals $(\alpha_h-1)^c(t)$, so that $(\alpha_h-1)^c=0$ on $X_*(H^0)$. Since $\alpha_h$ is unipotent and of finite order, we have $\alpha_h=1$.

Since, by Schur's lemma, the central torus $H^0$ acts on $V$ by scalars, the identity component of $H\cap\SL(V)$ is contained in the finite group $H^0\cap\SL(V)$ and $H\cap\SL(V)$ is finite. Since $\psi(s)\in[H,H]\subseteq H\cap\SL(V)$, the element $\psi(s)$ has finite order.
\end{proof}

\begin{proposition}\label{prop:nilpotent-mod-torsion for character varieties}
Let $\Gamma$ be a finitely presented group and let $\Gamma_{\mathrm{tors}}\trianglelefteq\Gamma$ be the subgroup generated by all elements of finite order. If $\Gamma/\Gamma_{\mathrm{tors}}$ is nilpotent, then the morphism
\[
    \mc X_G(\Gamma^\ab_\tf)_{0,\red}\to \mc X_G(\Gamma)_{0,\red}
\]
induced by the quotient map $\Gamma\to \Gamma^\ab_\tf$ is an isomorphism.

Consequently, the natural morphism
\[
    \mc X_T(\Gamma)_0/W\longrightarrow \mc X_G(\Gamma)_{0,\red}
\]
is finite and bijective, and it factors as an isomorphism onto the normalization of $\mc X_G(\Gamma)_{0,\red}$ followed by a homeomorphism.
\end{proposition}

\begin{proof}
Set $\overline\Gamma=\Gamma/\Gamma_{\mathrm{tors}}$. Since $\overline\Gamma$ is finitely generated and nilpotent, it is finitely presented. So by \Cref{lem:quotient by finite-order-in-rep-elements} applied to all elements of finite order, the natural morphism
\[
    \mc X_G(\overline \Gamma)_{0,\red}\to \mc X_G(\Gamma)_{0,\red}
\]
is an isomorphism. If $\gamma$ lies in the kernel of $\overline\Gamma\to \overline\Gamma^\ab_\tf$, then $\gamma^k\in[\overline\Gamma,\overline\Gamma]$ for some $k\geq1$, so that, by \Cref{lem:roots of unity}, $\psi(\gamma)^k$, and hence $\psi(\gamma)$, has finite order for every finite-dimensional irreducible representation $\psi$ of $\overline\Gamma$. Hence, we can apply \Cref{lem:quotient by finite-order-in-rep-elements} to all elements of the kernel of $\overline\Gamma\to \overline\Gamma^\ab_\tf$, thus obtain a natural isomorphism
\[
    \mc X_G(\overline \Gamma^\ab_\tf)_{0,\red}\to \mc X_G(\overline\Gamma)_{0,\red}.
\]
Noting that $\Gamma^\ab_\tf=\overline\Gamma^\ab_\tf$, which holds because the image of $\Gamma_{\mathrm{tors}}$ in $\Gamma^\ab$ is contained in $\mathrm{Tor}(\Gamma^\ab)$, shows that 
\(
    \mc X_G(\Gamma^\ab_\tf)_{0,\red}\longrightarrow \mc X_G(\Gamma)_{0,\red}
\)
is an isomorphism.

Since $\mc X_T(\Gamma)_0= \mc X_T(\Gamma^\ab_\tf)$, it follows that the natural morphism 
\(
    \mc X_T(\Gamma)_0/W\longrightarrow \mc X_G(\Gamma)_{0,\red}
\)
is equivalent to the natural morphism
\(
    \mc X_T(\Gamma^\ab_\tf)/W\longrightarrow \mc X_G(\Gamma^\ab_\tf)_{0,\red}
\)
As $\Gamma^\ab_\tf$ is free abelian, this morphism is finite and bijective by \cite[Proposition 7.1]{Sikora} and \cite[Theorem 3.4]{FLS25} (see also \cite{Thaddeus}). Because the algebraic torus $\mc X_T(\Gamma^\ab_\tf)$ is normal, the quotient $\mc X_T(\Gamma^\ab_\tf)/W$ is normal as well and the morphism factors as stated.
\end{proof}

\begin{proof}[Proof of \Cref{prop:reduction to T Betti}]
  Let $\Gamma=\pi_1(X)$. We have $M^T_\Betti(X)=\mc X_T(\Gamma)_{0}$ and $M^G_\Betti(X)=\mc X_G(\Gamma)_{0,\red}$ so that the assertion is a special case of \Cref{prop:nilpotent-mod-torsion for character varieties}.
\end{proof}

\begin{remark}\label{rem: nilpotent mod torsion}
Since Lemma \ref{lem:non-ab Hodge reduction} and the proofs of Theorems \ref{thm:P=W} and \ref{thm:geometric P=W} use the nilpotency of $\pi_1(X)$ only through Proposition \ref{prop:reduction to T Betti}, they hold line by line under the hypothesis of the proposition.
\end{remark}

\begin{remark}
    Imposing the much stronger conditions that $X$ is an abelian variety and $G$ is a classical group, Biswas--Florentino--Nozad  showed in \cite[Theorem 4]{BFN24} that  the morphism
    \[
    \Phi\colon M^{T}_\Betti(X)/W \to M^G_\Betti(X)
    \]
    is an algebraic isomorphism. 
    For $X$ an elliptic curve, descriptions of (the normalization of) $M^G_\Dol(X)$ of this kind appear in \cite{Thaddeus} and by Franco, Garc\'ia-Prada and Newstead \cite{FGN14, FGN19}.
\end{remark}

\begin{remark}
\label{rem:nilpotent examples}
Smooth projective varieties with nilpotent, non-abelian fundamental group exist: the first examples, of nilpotency class two, are due to Sommese and Van de Ven \cite{SommeseVandeVen86}, see also \cite{Campana95}. Proposition \ref{prop:nilpotent-mod-torsion for character varieties} refines, for the identity component of the character variety, the cohomological results of \cite[Corollary 1.4]{BergeronSilberman16} and \cite[Theorem 5.9]{FLS25}. If $G$ is a classical group, then $\mc X_G(\Z^r)_{0,\red}$ is normal by \cite{Sikora}, so that $M^G_\Betti(X)$, with its reduced structure, is normal for every $X$ with nilpotent fundamental group.
\end{remark}

\subsection{The Dolbeault moduli space }
\label{sec:preliminaries-dolbeault}
A $G$-Higgs bundle on a smooth projective variety $X$ is a pair $(E,\theta)$ consisting of a principal $G$-bundle $E$, together with a Higgs field $\theta\in H^0(\ad(E)\otimes_{\mc O_X} \Omega^1_X)$ satisfying the integrability condition $\theta\wedge \theta =0$. 
We denote by $M^G_\Dol(X)$ the connected component of the moduli space of semistable $G$-Higgs bundles on $X$ with vanishing characteristic classes  containing the trivial $G$-Higgs bundle.

The adjoint representation is trivial for $\G_m$, as a result a $\G_m$-Higgs bundle is a $\G_m$-bundle together with a global section of $\Omega^1_X$. In other words
\[
    M^{\G_m}_\Dol(X)= \Pic^0(X) \times H^0(X,\Omega^1_X)\ .
\]
It is easy to verify that
\[
    M^T_\Dol(X)= M^{\G_m}_\Dol(X) \otimes_\Z N \ .
\]
By extension of the structure group, one obtains a morphism $M^T_\Dol(X)\to M^G_\Dol(X)$ factoring through the quotient $M^T_\Dol(X)/W$. According to  \cite[Theorem 2.6]{Bi}, the induced morphism
\[
\Psi\colon    M^T_\Dol(X)/W \to M^G_\Dol(X)
\]
is injective.\footnote{Theorem 2.6 in \cite{Bi} is stated under the hypothesis that $X$ is an abelian variety. However, the argument for the injectivity of $\Psi$ works verbatim for arbitrary $X$ after  replacing the dual abelian variety by $\Pic^0(X)$.}

\begin{remark}
If $X$ is an abelian variety, it is shown in \cite[Theorem 2.6]{Bi} that $\Psi$ is  surjective as well. Moreover, \cite[Theorem 4]{BFN24} proves  $\Psi$  is an algebraic isomorphism if we assume in addition that $G$ is a classical group. Since we only need to work with much less restrictive hypotheses, we can go around the above assumptions  and show in Lemma \ref{lem:non-ab Hodge reduction} that $\Psi$ is a homeomorphism.  %through Proposition~\ref{prop:} and the non-abelian Hodge correspondence.     
\end{remark}
 
%%%%%%%%%%%%

\subsection{The correspondence}

The non-abelian Hodge correspondence states that there is a homeomorphism
\[
    \eta^G_X\colon M^{G}_\Betti(X)\to M^G_\Dol(X) \ .
\]
It turns out that, when $X$ has nilpotent fundamental group, the case of $\G_m$ fully determines what happens for an arbitrary reductive group $G$. 

\begin{lemma}\label{lem:non-ab Hodge reduction}
Let $X$ be an arbitrary smooth projective variety, $G$ a connected reductive algebraic group. Then we have a commutative diagram
\begin{equation}
\label{eq:eta for reductive group via tori}
\begin{tikzcd}
    M^{T}_\Betti(X)/W
    \arrow[r,"\eta_X^{T}/W"]
    \arrow[d,"\Phi"]&[2em]
        M^{T}_\Dol(X)/W
        \arrow[d,"\Psi"]\\
    M^G_\Betti(X)
    \arrow[r,"\eta_X^G"]&
        M^G_\Dol(X) \ .
\end{tikzcd}
\end{equation}
If $\pi_1(X)$ is nilpotent then the vertical arrows are homeomorphisms. Moreover, if $\widetilde M^G_\Dol(X)$ denotes the normalization of $M^G_\Dol(X)$, and $X$ has a nilpotent fundamental group, then the  homeomorphism $M^T_\Dol(X)/W\to M^G_\Dol(X)$ factors as the composite of an isomorphism $M^T_\Dol(X)/W\to \widetilde M^G_\Dol(X)$ and a homeomorphism $\widetilde M^G_\Dol(X)\to M^G_\Dol(X)$.
\end{lemma}
\begin{proof}
The diagram commutes by  the functoriality of the non-abelian Hodge correspondence \cite[\S6]{Simpson_nonabelianHodge}.
Next, assume that $\pi_1(X)$ is nilpotent. The vertical arrow on the left is a homeomorphism by Proposition \ref{prop:reduction to T Betti}. Since the horizontal arrows are also homeomorphisms, the vertical arrow on the right must be  one as well. In particular, it is a proper continuous map. By GAGA the morphism $M^T_\Dol(X)/W\to M^G_\Dol(X)$ is also proper as a morphism of schemes, and by Zariski Main Theorem it is finite. The space $M^T_\Dol(X)\cong (\Pic^0(X)\times H^0(X,\Omega^1_X))\otimes_\Z N$ is normal, implying that $M^T_\Dol(X)/W$ is normal as well. It now follows that the induced morphism $M^T_\Dol(X)/W\to \widetilde M^G_\Dol(X)$ is an isomorphism.
\end{proof}

%%%%%%%%%%%%%%%%%%%%%%%%%%%%%%%%%%%%%%%%%%%%%%%%%%%
%% SECTION 3
%%%%%%%%%%%%%%%%%%%%%%%%%%%%%%%%%%%%%%%%%%%%%%%%%%%
\section{The spectral data morphism and the cohomological $P=W$ conjecture}

We briefly recall and unwind the notion of the spectral data morphism introduced by  Chen--Ng\^o \cite{CN20}.  Let $\mathfrak{t}$ and $\mathfrak{g}$ denote the Lie algebras of $T$ and $G$. We denote by $\mathfrak C_G^d$ the commuting scheme of $d$ elements in $\mathfrak g$ with pairwise vanishing Lie brackets, that is, the scheme theoretic fiber  of the commutator map 
\[
    \mathfrak{g}^d\rightarrow \prod_{1\leq i<j\leq d}\mathfrak{g},\quad (\theta_1,\ldots,\theta_d)\mapsto ([\theta_i,\theta_j])_{i<j}
\]
over $0$.  There is an action of $\GL_d\times G$ on $\mathfrak g^d=\A^d\otimes_k \mathfrak g$, where $G$ acts on $\mathfrak g$ with the adjoint representation and $\GL_d$ on $\A^d$ via the vector representation. The commuting scheme $\mathfrak C_G^d$ is invariant under this action, hence carries an induced action. 

The stack quotient $[\mathfrak C_G^d/\GL_d\times G]$ parametrizes triples $(\mc V, \mc E, \theta)$ consisting of a vector bundle $\mc V$ of rank $d$, a principal $G$-bundle $\mc E$, and a morphism of vector bundles $ \theta\colon\mc V^\vee\rightarrow\ad(\mc E)$ satisfying the integrability condition $[\theta(\xi),\theta(\eta)]=0$ for all local sections $\xi,\eta$ of $ \mathcal{V}^{\vee}$.

To obtain $G$-Higgs bundles on $X$, we must set $d=\dim X$ and take $\mc V=\Omega^1_X$. In this case, the map $\theta$ is a section of $\ad(\mathcal{E})\otimes \Omega^1_X$, and, writing $\theta=\sum_i\theta_idz_i$ in local coordinates we have $\theta\wedge\theta=\sum_{i<j}[\theta_i,\theta_j]dz_i\wedge dz_j$. Therefore,  the integrability condition above is exactly the condition $\theta\wedge \theta=0$. To force $\mathcal{V}=\Omega_X^1$, we must take the fiber product
\[
    \mc H_X\coloneqq X\times_{\mbb B\GL_d}[\mathfrak C_G^d/\GL_d\times G]
\]
of the morphism $X\xrightarrow{\Omega^1_X} \mbb B\GL_d$ defined by the cotangent sheaf $\Omega_X^1$ and the map $[\mathfrak C_G^d/\GL_d\times G]\to \mbb B \GL_d$ that forgets everything but the vector bundle. We view $\mc H_X$ as a stack over $X$; by construction its sections over $X$ are precisely the $G$-Higgs bundles on $X$. 

The moduli stack $\mc M^G_\Dol(X)$ of $G$-Higgs bundles on $X$ is then the Weil restriction $\Res_{X/k}(\mc H_X)$ of $\mc H_X$ to $\Spec k$. In other words,  for a test $k$-scheme $S$, the $S$-valued points of $\mc M^G_\Dol(X)$  are  the morphisms $X\times_k S \to \mathcal{H}_X$, which are precisely the families of $G$-Higgs bundles over $S$. 
The moduli space $M_{\Dol}^G(X)$ is the identity component of the good moduli space of the open substack of $\mc M^G_\Dol(X)$ consisting of semistable $G$-Higgs bundles \cite{SimI}.

The spectral data morphism defined by Chen and Ng\^o \cite{CN20} removes the principal $G$-bundle from a $G$-Higgs bundle in the minimal way by replacing the stack theoretic quotient by $G$ by a categorical quotient. Formally, they take as the Hitchin base $B^G_X$ the space
\begin{equation}
\label{eq:def of Hitchin base}
    B^G_X\coloneqq \Res_{X/k}(X\times_{\mbb B\GL_d} [(\mathfrak C^d_G/\!\!/G)/\GL_d])\ ;
\end{equation}
namely, the space of sections of the fiber bundle on $X$ with fiber $\mathfrak C^d_G/\!\!/G$ obtained by twisting with the frame bundle of $\Omega_X^1$.
Choosing a $\GL_d$-equivariant closed embedding of the affine scheme $\mathfrak C^d_G/\!\!/G$ into a finite dimensional representation $V$ of $\GL_d$, we see that $B^G_X$ is a closed subscheme of the finite dimensional affine space of global sections of the vector bundle on $X$ associated with $V$; in particular, $B^G_X$ is an affine scheme of finite type over $k$. The quotient morphism $[\mathfrak C^d_G/\GL_d\times G]\to[(\mathfrak C^d_G/\!\!/G)/\GL_d]$ induces a morphism $\mc M^G_\Dol(X)\to B^G_X$. Since good moduli spaces are universal for morphisms to algebraic spaces \cite{Alper13}, its restriction to the open substack of semistable, topologically trivial $G$-Higgs bundles descends to a morphism
\[
\sd^G_X\colon M^G_\Dol(X)\longrightarrow B^G_X\ ,
\]
which Chen--Ng\^o call the \emph{spectral data morphism}.

The canonical inclusion $\mathfrak{t}\subseteq \mathfrak{g}$ induces a $W$-equivariant map $\mathfrak{t}^d=\mathfrak{C}^d_T\to \mathfrak{C}^d_G$ and, in turn, a morphism
\[
    \mathfrak t^d/W \to \mathfrak C^d_G/\!\!/G\ .
\]
Note that the above map is an isomorphism modulo  a conjectural generalization of Chevalley restriction theorem. The latter holds  true for $d=1$, where it is precisely Chevalley's restriction theorem.
For $d\geq 2$, the result is known for $G=\GL_n(\C)$  \cite{CN20} and $G=\Sp_{2n}(\C)$ \cite{CN21}, extending earlier results of Hunziker \cite{H97} and Joseph \cite{J97}; see also \cite{Xia24} for classical groups in positive characteristic. None of the results below depends on the generalized Chevalley restriction theorem, with the exception of Lemma \ref{lem:trivial cotangent}, where it appears as a hypothesis.

\begin{remark}
If one chooses homogeneous generators $c_1,\ldots, c_h$ of degrees $e_1,\ldots, e_h$, respectively,  on the affine space $\mathrm{Spec} (k[\mathfrak{g}]^G)=\mathfrak g/\!\!/G\cong \mathfrak t/W=\mathrm{Spec}(k[\mathfrak{t}]^W)$,  then the Weyl polarizations 
\[
    \mathfrak g^d(k)\ni(\theta_1,\ldots, \theta_d)\mapsto c_i(x_1\theta_1+\cdots+x_d \theta_d)\in k[x_1,\ldots, x_d]_{e_i}
\]
define a $\GL_d$-equivariant  and $G$-invariant polarization morphism
\[
    \mathrm{pol}\,\colon\, \mathfrak C^d_G\to \prod_{i=1}^h S^{e_i}\A^d\ ,
\]
where $S^{e_i}\A^d$ is the $e_i$-th symmetric power of $\A^d$, that is, the vector space $k[x_1,\ldots, x_d]_{e_i}$ of degree $e_i$ polynomials in $d$ variables viewed as an affine space.

Note that the morphism $\mathrm{pol}$ factors through $\mathfrak{C}^d_G/\!\!/G$, thus it induces a morphism 
\[
    \iota_X^G\colon B^G_X\rightarrow 
        \mathrm{Res}_{X/k}\left( X\times_{\mathbb{B}\GL_d}\left[\prod_{i=1}^h S^{e_i}\A^d \Big/\GL_d \right]\right)\cong
            \prod_{i=1}^h H^0(S^{e_i}\Omega^1_X) 
                \eqqcolon\mathcal{A}^G_X\ .
\]
The composition $h_X^G\coloneqq\iota_X^G\circ \sd_X^G$ is the Hitchin morphism, which sends a Higgs pair $(E,\theta)$ to $(c_1(\theta),\ldots, c_h(\theta))$.
The Hitchin map is proper \cite{SimI} and, since $\iota_X^G$ is separated, $\sd_X^G$ is proper.
However, in general $h_X^G$ fails to be surjective for dimensional reasons. The spectral data morphism $\sd$ may be thought of as the Hitchin morphism with its target replaced by a space of spectral data that takes into account the integrability condition.

In \cite{CN20}, the Hitchin base is defined via  replacing $\mathfrak{C}^d_G$ by its scheme theoretic image $\mathcal{B}^G_X$ under $\mathrm{pol}$ in \eqref{eq:def of Hitchin base}, so that the corresponding map $\mathcal{B}^G_X\hookrightarrow \mathcal{A}^G_X$ is in fact a closed immersion. While this makes the connection to the Hitchin map more transparent, for our purposes the difference between the two definitions is immaterial. In fact, the relevant maps between the various bases are finite and injective (cf.\ Lemma \ref{lem:finiteinjectivebetaxg}) and, as a result, they give rise to the same perverse filtration on the cohomology of $M_{\Dol}^G(X)$ (cf.\ Proposition \ref{prop:allperverseequal}).
\end{remark}

In case $G=T$ is a maximal torus, we can make the description of spectral data morphisms explicit.
\begin{lemma}
\label{lem:functoriality of spectral data}
There is a natural identification $B_X^T\cong H^0(X,\Omega^1_X)\otimes_\Z N$ such that the spectral data morphism $\sd^T_X$ becomes  the projection
\[
    M^T_\Dol(X)\cong (\Pic^0(X)\oplus H^0(X,\Omega^1_X))\otimes_\Z N\to H^0(X,\Omega^1_X)\otimes_\Z N\ .
\]
Moreover, the diagram
\[
\begin{tikzcd}
    M_\Dol^T(X) \arrow[r,"\sd_X^T"]\arrow[d]& 
        B^T_X  \arrow[d]\\
    M_\Dol^G(X) \arrow[r,"\sd_X^G"]&
        B^G_X
\end{tikzcd}
\]
commutes.
\end{lemma}
\begin{proof}
The functoriality of the spectral data morphism in the reductive group, and hence the commutativity of the diagram, is immediate from the definition. In the torus case, we have that $T$ acts trivially on $\mathfrak C^d_T=\mathfrak t^d$ and thus that $\mathfrak C^d_T/\!\!/T=\mathfrak t^d=\mathfrak t\otimes_k \A^d$. It follows that for any $k$-scheme $S$, the $X\times_k S$-valued points of  
\[
    X\times_{\mbb B\GL_d} [(\mathfrak C_T^d/\!\!/T)/\GL_d]
\]
are in natural bijection with 
\[
    H^0(X\times_k S,pr_X^*\Omega^1_X\otimes_k\mathfrak t) = (H^0(\Omega^1_X)\otimes_k \mathfrak t)\otimes_k H^0(S,\mc O_S) =(H^0(\Omega^1_X)\otimes_\Z N)\otimes_k H^0(S,\mc O_S)\ ,
\]
which, in turn,  agrees with $(H^0(\Omega^1_X)\otimes_\Z N)(S)$.  
To apply the spectral data morphism for $T$-Higgs bundles, we replace the stack quotient by $T$ through  the categorical quotient by $T$. As $T$ acts trivially on $\mathfrak{C}^d_T=\mathfrak t^d$, this amounts to forgetting the $T$-bundle.
\end{proof}

The Weyl group action on $\mathfrak t$ induces an action on $B^T_X$. In the identification of $B^T_X$ with $H^0(X,\Omega^1_X)\otimes_\mathbb{Z} N$, this coincides with  the action induced by the $W$-action on the cocharacter lattice $N$.  Since the morphism $B^T_X\rightarrow B_X^G$ is $W$-invariant, it  descends to a morphism
\[
    \beta^G_X\colon \widetilde B^G_X\coloneqq B^T_X/W \to B^G_X\ .
\]

\begin{remark} In general, $\beta^G_X$ is far from being an isomorphism. One can take for example  a curve $X$ of genus $g\geq 2$, and $G=\SL_2$. In this case 
\[
    \mathcal B^G_X=\mathcal{A}^G_X\cong H^0(X,K_X^2)\cong \mathbb{A}^{3g-3}\ ,
\]
while $\widetilde{B}^G_X\cong H^0(X,K_X)/\{\pm 1\}$ has dimension $g$.
\end{remark}

For abelian varieties, the following lemma shows that $\beta_X^G$ is an isomorphism provided the generalized Chevalley  restriction theorem holds. 
 
\begin{lemma}\label{lem:trivial cotangent}
Assume that $\Omega^1_X$ is trivial, and the generalized Chevalley restriction theorem holds. Then $\widetilde B^G_X\to B^G_X$ is an isomorphism.    
\end{lemma}
\begin{proof}
The fiber product $X\times_{\mbb B\GL_d} [(\mathfrak C_G^d/\!\!/G)/\GL_d]$ is a fiber bundle over $X$ with fiber $\mathfrak C_G^d/\!\!/G$ and transition morphisms induced by those of $\Omega^1_X$. So, if $\Omega^1_X$ is trivial, and the generalized Chevalley restriction theorem holds, then  we have
\[
    B^G_X  \equ  \Res_{X/k}(X\times_k (\mathfrak t^d/W))= \mathfrak t^d/W\ .
\]
Applying this in the case $G=T$ we obtain $B^T_X= \mathfrak t^d$. The morphism 
\[
    \mathfrak t^d=B^T_X\to B^G_X=\mathfrak t^d/W
\]
is precisely the quotient morphism.
\end{proof}

For the purposes of our paper, however, we do not need $\beta_X^G$ to be an isomorphism. It suffices if  $\beta_X^G$ is finite and injective, and this property holds without any further assumption on $X$ and $G$.

\begin{lemma}\label{lem:finiteinjectivebetaxg}
The composition $\iota_X^G\circ \beta_X^G\colon\widetilde B_X^G\rightarrow \mathcal{A}_X^G$ is finite and injective on $k$-points. As a consequence,   $\beta_X^G$ is finite and injective on $k$-points as well.
\end{lemma}

\begin{proof}
For a $k$-point $p=(x,v)$ of the total space $\operatorname{Tot}(T_X)$ of the tangent bundle  (that is, $x\in X$ and $v\in T_xX$), and for $\alpha\in B^T_X=H^0(X,\Omega^1_X)\otimes_k\mathfrak t$, we write $\alpha(p)\coloneqq\langle\alpha(x),v\rangle\in\mathfrak t$. Similarly, for $s\in H^0(X,\Sym^e\Omega^1_X)$ we write $s(p)\in k$ for the value at $v$ of the degree $e$ polynomial $s(x)$ on $T_xX$. By construction, the composite
\[
    c\colon B^T_X\longrightarrow\widetilde B^G_X\xrightarrow{\ \iota^G_X\circ\beta^G_X\ }\mc A^G_X
\]
sends $\alpha$ to $(c_1(\alpha),\ldots,c_h(\alpha))$, where $c_i(\alpha)(p)=c_i(\alpha(p))$ for every $p$.
 
We first prove the  injectivity of the composition. Let $\alpha,\alpha'\in B^T_X$ with $c(\alpha)=c(\alpha')$. For every $p$, the elements $\alpha(p)$ and $\alpha'(p)$ of $\mathfrak t$ then have the same image under the quotient map $\mathfrak t\to\mathfrak t/W$, as $c_1,\ldots,c_h$ generate $k[\mathfrak t]^W$; hence they lie in the same $W$-orbit. It follows that the closed subsets
\[
    Z_w\coloneqq\{p\in\operatorname{Tot}(T_X):\alpha'(p)=w\alpha(p)\}\ ,\qquad w\in W\ ,
\]
cover all $k$-points of $\operatorname{Tot}(T_X)$, and hence $\operatorname{Tot}(T_X)$ itself. Since $X$ is irreducible, so is $\operatorname{Tot}(T_X)$, and therefore $Z_w=\operatorname{Tot}(T_X)$ for some $w\in W$, that is, $\alpha'=w\alpha$. As $k$ is algebraically closed, the $k$-points of $\widetilde B^G_X=B^T_X/W$ are the $W$-orbits in $B^T_X$, so $\iota^G_X\circ\beta^G_X$ is injective on $k$-points.
 
Next, we  prove finiteness. A non-zero $\alpha\in B^T_X$ satisfies $\alpha(p)\neq0$ for some $p$, and $B^T_X$ is finite-dimensional. Hence we can choose finitely many points $p_1,\ldots,p_N$ such that the linear map
\[
    \mathcal{E}\colon B^T_X\longrightarrow\mathfrak t^N\ ,\quad\alpha\longmapsto\big(\alpha(p_1),\ldots,\alpha(p_N)\big)
\]
is injective. Consider the linear map $\varepsilon\colon\mc A^G_X\to\A^{hN}$, $(s_i)_i\mapsto(s_i(p_j))_{i,j}$, and the morphism $\chi=(c_1,\ldots,c_h)\colon\mathfrak t\to\A^h$. Then
\[
    \varepsilon \circ c=\chi^N\circ \mathcal E\ .
\]
The morphism $\chi$ is finite, since $k[\mathfrak t]$ is integral over $k[\mathfrak t]^W=k[c_1,\ldots,c_h]$, $W$ being finite. As $\mathcal E$ is a closed immersion, $\varepsilon\circ c$ is finite, and since $\varepsilon$ is separated, so is $c$.
%\footnote{Recall that if $f\colon U\to V$ and $g\colon V\to V'$ are morphisms such that $g\circ f$ is finite and $g$ is separated, then $f$ is finite: indeed, $f$ factors as the graph $U\to U\times_{V'}V$, which is a closed immersion, followed by the base change of $g\circ f$ along $g$.} 
Thus $k[B^T_X]$ is a finite $k[\mc A^G_X]$-module, and so is its submodule $k[\widetilde B^G_X]=k[B^T_X]^W$, because $k[\mc A^G_X]$ is noetherian. In other words, $\iota^G_X\circ\beta^G_X$ is finite. Finally, since $\iota^G_X$ is separated, $\beta^G_X$ is finite as well.
\end{proof}

\subsection{Spectral data morphism and the  associated perverse filtration for varieties with nilpotent fundamental group}

Our treatment of the spectral data morphism has not put  any assumption on $X$ so far. We now specialize to the case of interest, where $\pi_1(X)$ is a nilpotent group. By definition, the spectral data morphism $\sd^T_X$ is $W$-equivariant, and thus induces a morphism 
\[
    \widetilde\sd_X^G\colon M^T_\Dol(X)/W\cong \widetilde M_{\Dol}^G(X) \rightarrow \widetilde B^G_X=B^T_X/W 
\]
on the quotients by $W$. This morphism is surjective and proper as $\sd^T_X$ is a projection with compact fiber $\Pic^0(X)\otimes_{\mathbb{Z}}N$ by Lemma \ref{lem:functoriality of spectral data}. 
We obtain a commutative diagram 
\begin{equation}
	\label{eq:spectral data diagram}
	\begin{tikzcd}
		\widetilde M^G_\Dol(X)\cong M^T_\Dol(X)/W
		\arrow[r,"\widetilde\sd^G_X"]
		\arrow[d,"\cong_{homeo}"]&
		\widetilde B^G_X
		\arrow[d,"\beta^G_X"]&\\
		M^G_\Dol(X)
		\arrow[r,"\sd^G_X"]&
		B^G_X
		\arrow[r,"\iota^G_X"]&
		\mc A^G_X\ ,
	\end{tikzcd}
\end{equation}
where the vertical arrow on the left is the normalization morphism and the composite of the bottom arrows is the Hitchin morphism $h_X^G$. 
Since the normalization map $\nu\colon \widetilde M_{\Dol}^G(X) \to M_{\Dol}^G(X)$ is a homeomorphism by Lemma \ref{lem:non-ab Hodge reduction}, we can view $\widetilde \sd^G_X$ as a map 
\[
    M_{\Dol}^G(X) \rightarrow \widetilde B^G_X\ ,
\]
which, with a slight abuse of notation, we still denote by $\widetilde \sd^G_X$.

The morphism $\widetilde \sd^G_X$ is the \textquotedblleft minimal\textquotedblright\ map that is compatible with the extension of the structure group from $T$ to $G$, and it is the one we will work with in the proof of Theorem A.
However, as we now show, we are not losing relevant information by making this choice. The perverse filtrations associated with the morphisms $\sd_X^G,\widetilde \sd_X^G,$ and $h^G_X$ all coincide.

\begin{proposition}\label{prop:allperverseequal}
Suppose that $\pi_1(X)$ is nilpotent. Then the morphisms $\sd_X^G,\widetilde \sd_X^G,$ and $h^G_X$ are proper, and they induce the same perverse filtration on $H^*(M_{\Dol}^G(X),\mathbb{Q})$. More precisely, if we identify $H^*(M_{\Dol}^G(X))$ with $H^*(\widetilde M_{\Dol}^G(X))$ via $\nu^*$, we have 
\[
    P_p^{h^G_X}H^*(M_{\Dol}^G(X),\mathbb{Q})=P_p^{\sd^G_X}H^*(M_{\Dol}^G(X),\mathbb{Q})=
    (\nu^*)^{-1}P_p^{\widetilde\sd^G_X}H^*(\widetilde M_{\Dol}^G(X),\mathbb{Q})\ .
\]
\end{proposition}
\begin{proof}
The proof follows directly from Lemma \ref{lem:finiteinjectivebetaxg} and Lemma \ref{lem:perversecompostion} applied with $f=\widetilde{\sd}^G_X$, $g=\beta_X^G$ and $\nu$ given by the normalization. 
\end{proof}

\begin{remark}[Convention regarding perverse filtrations]
    Because of Proposition \ref{prop:allperverseequal}, we will omit the superscript in the definition of the perverse filtration for the rest of the paper,  and simply write $P_pH^*(M_{\Dol}^G(X),\mathbb{Q})$ for any of the three. 
\end{remark}

\subsection{Cohomological $P=W$ for varieties with nilpotent fundamental group}

In  this subsection we prove Theorem \ref{thm:P=W}, our  cohomological $P=W$ result over varieties with nilpotent fundamental group. 
The latter condition allows us to reduce cohomological $P=W$ from $G$ to $T$. In the particular situation  where $G=T$ is a torus, we need not assume that $\pi_1(X)$ is nilpotent, and cohomological $P=W$ holds over an arbitrary smooth projective variety. 

\begin{proposition}[$P=W$ for $G =T$ and arbitrary $X$]
\label{prop:Pw_rank1_arbitrary}
Let $X$ be any smooth complex projective variety, and let $G = T$ be an algebraic torus of rank $n$. Then the cohomological $P=W$ identity holds unconditionally for $X$:
\[
P_k H^d(M_{\Dol}^{T}(X), \mathbb{Q}) = W_{2k} H^d(M_{\Betti}^{T}(X), \mathbb{Q}) = W_{2k+1} H^d(M_{\Betti}^{T}(X), \mathbb{Q}) \quad \forall k, d \ge 0\ .
\]
\end{proposition}

\begin{proof}
By Lemma \ref{lem:functoriality of spectral data}, the spectral data morphism
\[
    \sd^T_X\colon M^T_\Dol(X)\cong (\Pic^0(X)\oplus H^0(X,\Omega^1_X))\otimes_\Z N \to H^0(X,\Omega^1_X)\otimes_\Z N \cong B^T_X
\]
is the projection to the second factor.
By the K\"unneth formula, the pushforward $R(\sd^{T}_X)_* \mathbb{Q}_{M_{\Dol}^{T}(X)}$ decomposes into a direct sum of shifted constant sheaves 
\[
    R(\sd^{T}_{X})_* \mathbb{Q}_{M_{\Dol}^{T}(X)}\cong \bigoplus_{i=0}^{2gn} H^i(\operatorname{Pic}^0(X)\otimes_\Z N, \mathbb{Q}) \otimes \mathbb{Q}_{B^T_X}[-i]\ ,
\]
where $2g=b_1(X)$ is the first Betti number of $X$.\\
Consequently, the perverse filtration on $H^*(M_{\Dol}^{T}(X), \mathbb{Q})$ induced by $\sd^T_{X}$ is governed purely by the cohomological degree on $\operatorname{Pic}^0(X)\otimes_{Z}N$.

On the Betti side, any representation $\rho\colon \pi_1(X) \to T$ factors through the group $H_1(X, \mathbb{Z})$, yielding $M_{\operatorname{Betti}}^{T}(X)=\mathrm{Hom}(H_1(X, \mathbb{Z}),T)_0$, which is an algebraic torus. By Deligne's mixed Hodge theory for algebraic tori, the cohomology group $H^i(M_{\Betti}^{T}(X), \mathbb{Q})$
%\cong \bigwedge^i H^1(\, \mathbb{Q})$ 
is pure of weight $2i$. 

It follows that both the perverse and the weight filtration are determined fully by cohomological degree; matching the degrees proves the assertion.
\end{proof}
 
Having collected the necessary auxiliary results, we can now move on to proving the first main result of the paper, that is that the weight filtration on $M^G_\Betti(X)$ can be matched with the perverse filtration induced by $\widetilde\sd^G_X$ when $X$ has nilpotent fundamental group. 

\begin{proof}[Proof of Theorem \ref{thm:P=W}]
%By Lemma~\ref{lem:reduction to pi1 free abelian} we can and will assume that $\pi_1(X)$ is a free abelian group.  
We reduce the statement to the case $G=T$. To this end,  consider the commutative diagram
\[
\begin{tikzcd}
    M_\Betti^{T}(X)
    \arrow[d]
    \arrow[r,"\eta_X^{T}"]& [3em]
        M_\Dol^{T}(X)
        \arrow[d]
        \arrow[r,"\sd_X^{T}"]&[3em]
            B^T_X
            \arrow[d]\\
    M_\Betti^{T}(X)/W 
    \arrow[r,"\eta_X^{T}/W"]
    \arrow[d] &
        M_\Dol^{T}(X)/W
        \arrow[d]
        \arrow[r,"\sd^{T}_X/W"]&
            \widetilde B^G_X
            \arrow[d,"="]\\
    M_\Betti^G(X)
    \arrow[r,"\eta_X^G"]&
        M_\Dol^G(X)
        \arrow[r,"\widetilde\sd^G_X"]&
            \widetilde B^G_X\ .
\end{tikzcd}
\]
The horizontal arrows on the left are all homeomorphisms by the non-abelian Hodge correspondence, and thus induce isomorphisms on cohomology of the left and  middle spaces. We claim that these isomorphisms identify the weight filtrations of the spaces on the left with the perverse filtrations of the spaces in the middle column, with respect to the maps to the right column. Once we show this claim for the last row, we have proved the assertion.

By Proposition \ref{prop:Pw_rank1_arbitrary},  the statement holds for the first row. The morphisms in the first row are $W$-equivariant and the second row is obtained from the first row taking the quotient by $W$. Therefore, the claim for the second row follows from Proposition \ref{prop:quotient-filtration} and the functoriality of weight filtration under finite group actions. Finally, the lower vertical arrows are all algebraic and homeomorphisms by Lemma \ref{lem:non-ab Hodge reduction}. By Lemma \ref{lem:perversecompostion}(ii) it follows that the perverse filtrations of the middle spaces are identified by the middle vertical arrow. Because morphisms of mixed Hodge structures are strictly compatible on filtrations \cite[Theorem 2.3.5]{Del71}, it also follows that the left lower vertical arrow induces an isomorphism on weight filtrations (cf.\cite[Theorem 5.8]{FLS25}). Therefore, the claim for the second row implies the claim for the third row.
\end{proof}

Next, recall the Lefschetz hyperplane theorem for fundamental groups \cite[Theorem 3.1.21]{PAG_I}. % original source?

\begin{theorem}[Lefschetz]\label{thm:Lefschetz}
Let $X$ be a smooth projective variety, $D$ an effective ample divisor on $X$. If $\dim X\geq 3$ then  inclusion induces an isomorphism of groups 
\[
\pi_1(D) \to \pi_1(X)\ . 
\]
\end{theorem}

\begin{example}\label{eg:subvarieties of general type}
 In particular, if $X$ is a variety of dimension at least three with a nilpotent fundamental group, and $|A|$ a very ample linear series on $X$, then a general element $D\in |A|$ will also be a smooth projective variety having  an nilpotent fundamental group.  As in \cite[Example 3.1.31]{PAG_I},  the phenomenon extends by induction to subvarieties  $V\subseteq X$ that are general complete intersections of very ample divisors on $X$.  This gives a wide class of examples where Theorem \ref{thm:P=W} applies. 

If $X$ is in fact an abelian variety then complete intersections of ample hypersurfaces will be of general type. 
\end{example}

\begin{corollary}\label{cor:general type}
   For every connected complex reductive group $G$, there exist varieties of general type   in every dimension $n\geq 2$ such that 
   the  non-abelian Hodge correspondence 
\[
    \eta_X^G\colon M^G_\Betti(X)\to M^G_\Dol(X)
\]
induces an isomorphism
\[
    P_k H^*(M^G_\Dol(X),\Q) \xrightarrow{\cong}
    W_{2k} H^*(M^G_\Betti(X),\Q) \,=\,
    W_{2k+1} H^*(M^G_\Betti(X),\Q)\ ,
\]
where we consider  the perverse filtration  with respect to the spectral data morphism.
\end{corollary}
 
\begin{remark}
As opposed to all examples  to give rise to $P=W$ phenomena on the appropriate moduli spaces that were known before our current work, there are plenty of varieties that have nilpotent fundamental groups but are not aspherical. 
\end{remark}

\begin{example}[Non-nilpotent examples]\label{ex:non-nilpotent}
Let $H$ be a finite group that is not nilpotent. By Godeaux--Serre there is a smooth projective variety $Y$ with $\pi_1(Y)\cong H$. Let $Z$ be an abelian variety, or more generally any smooth projective variety with nilpotent fundamental group. Then $X=Y\times Z$ has fundamental group $H\times\pi_1(Z)$, which is not nilpotent, and is infinite as soon as $b_1(Z)>0$. Its quotient by the subgroup generated by elements of finite order is a quotient of $\pi_1(Z)$, hence nilpotent, so by Proposition \ref{prop:reduction to T Betti} and Remark \ref{rem: nilpotent mod torsion} the conclusions of Theorems \ref{thm:P=W} and \ref{thm:geometric P=W} hold for $X$.
\end{example}

%%%%%%%%%%%%%%%%%%%%%%%%%%%%%%%%%%%%%%%%%%%%%%%%%%%
%% SECTION 4
%%%%%%%%%%%%%%%%%%%%%%%%%%%%%%%%%%%%%%%%%%%%%%%%%%%

\section{The geometric $P=W$ conjecture}
\label{sec:geometric P=W}
The current section is devoted to the proof of a  geometric $P=W$ result over smooth projective varieties with nilpotent fundamental group. The argument goes through a similar statement for toric varieties, which was motivated by \cite{MMS_geometricP=W}.  

\subsection{The construction of the  evaluation map}
\label{sec:evaluation map}
Let $X$ be a smooth variety, and let $\overline X$ be a smooth compactification of $X$ with simple normal crossings (snc) boundary $\partial X=\sum_{i\in I}D_i$. We make the following assumption:
\begin{equation}
\label{eq:irrecubile condition}
     \text{\emph{All non-empty intersections of components of $\partial X$ are irreducible.}}
\end{equation}
    
Let $\Delta(\partial \overline X)$ denote the dual complex of $X$ in the sense of \cite{dFKX17}. The homotopy type of $\Delta(\partial \overline X)$ only depends on $X$ \cite{dFKX17}, so we may also denote it by $\Delta_X$.  For each component $D_i$ of $\partial X$, we can choose a \emph{rug function}, that is, a proper, real semialgebraic map
\[
    \beta^X_i\colon \overline X\to \R_{\geq 0}
\]
with $(\beta^X_i)^{-1}\{0\}=D_i$. One can associate a \emph{semialgebraic neighborhood} $N_i(X)$ to $\beta_i^X$ by a choice of $\delta>0$,  and by setting 
\[
N_i(X) \coloneqq (\beta^X_i)^{-1}[0, \delta)\ .
\]
We suppress $\delta$ from the notation for simplicity. Taking $\delta$ small enough, we may and will  make the following assumption:
\begin{equation}
\label{eq:condition dual complex}
 \text{\emph{For all }} J\subseteq I\text{\emph{, we have }}\bigcap_{i\in J} N_i(X) \neq \emptyset \text{\emph{ if and only if }} \bigcap_{i\in J} D_i\neq \emptyset\ .
\end{equation}
Semialgebraic neighborhoods are in fact unique up to isotopy \cite[Proposition 3.5]{Durfee_neighborhoods}.
Taking a union over all $N_i(X)$ and removing the boundary, we obtain a \emph{neighborhood at infinity }
\[
    N^*(X)= \bigcup_{i\in I} N_i(X) \setminus \partial X \ , 
\]
which is covered by the sets $N^*_i(X)= N_i(X)\setminus \partial X$.

If $\mathfrak f\colon \R_{\geq 0}\to [0,1]$ denotes the cutoff function that is $1$ at $0$, $0$ on $[\delta,\infty)$, and linear on $[0,\delta]$, we obtain functions 
\begin{equation}
\label{eq:partition of unity}
    \varphi^X_i\coloneqq \frac{\mathfrak f\circ \beta^X_i}{\sum_{j\in I} \mathfrak f\circ \beta^X_j}\ ,
\end{equation}
    
which form a partition of unity subordinate to the cover $\{N^*_i(X)\}_{i\in I}$ of $N^*(X)$. Combining the $\varphi^X_i$, we obtain the \emph{evaluation map}
\[
    \ev^X\colon N^*(X)\to \R^I,\; x\mapsto (\varphi^X_i(x))_{i\in I}\ .
\]
By assumptions \labelcref{eq:irrecubile condition} and \labelcref{eq:condition dual complex}, the image of $\ev^X$ will be contained in the standard embedding of $\Delta(\partial X)$ in $\R^I$. Moreover, we may in fact choose the rug functions in a way that the partition of unity is \emph{full}, that is, such that $\ev^X$ surjects onto $\Delta(\partial X)$. The evaluation map $\ev^X\colon N^*(X)\to \Delta(\partial X)$ clearly depends on a variety of choices, but its homotopy type does not, at least not if $\Delta(\partial X)$ can be chosen to be a homotopy sphere. Moreover, we need not assume that $X$ is smooth and that $\overline X$ is an snc model of $X$: it suffices for $\overline X$ to be a dlt model \cite{dFKX17}.

If a finite group $H$ acts on $X$ and this action extends to $\overline X$, there is an induced action on $I$ given by $D_{hi}=hD_i$,  and hence, on $\Delta(\partial X)$. If we replace each rug function $\beta^X_i$ by the $H$-equivariant functions
\begin{equation*}
\label{eq:sum of shifts}
    \widetilde \beta^X_i(x)= \sum_{h\in H} \beta^X_{hi}(hx)\ , 
\end{equation*}
then the set-theoretic neighborhood at infinity $N^*(X)$ is $H$-invariant and the evaluation map $\ev^X$ is $H$-equivariant.

\subsection{Geometric $P=W$ for toric  varieties}
\label{sec:evaluation map polytopes}
We elaborate on the observation that the construction of Section \ref{sec:evaluation map} also works when $X=\mathring{P}$ is the interior of a full-dimensional polytope $\overline X=P$ in $K_\R$ for some lattice $K$. We may assume  without loss of generality  that $0\in \mathring P$. The boundary 
\[
\partial P \equ \bigcup_{i\in I} D_i
\]
is a union of the facets $D_i$ of the polytope $P$, and for $\beta^P_i$ we can take the absolute value $|l_i|$ of the unique affine function $l_i$ on $K_\R$ which vanishes precisely on the affine span of $D_i$ and satisfies $l_i(0)=1$.
For the dual complex $\Delta(\partial P)$ one has to be careful, however, as the analogous definition only works well if $P$ is simple, and does not appropriately capture the face lattice of $P$ otherwise. Instead, we define $\Delta(P)\coloneqq \partial P^\Delta$, where
\[
P^\Delta\coloneqq \{m\in K^\vee_\R : m|_P\leq 1\}
\]
is the polar polytope of $P$ in the dual space associated to the dual lattice $K^\vee=\Hom(K,\Z)$. The polar has the property that it has a vertex $v_i^\Delta\in P^\Delta$ for every facet $D_i$, $i\in I$, of $P$ and for $J\subset I$, the vertices $\{v_j^\Delta\}_{j\in J}$ are contained in a facet if and only if $\bigcap_{j\in J}D_j\neq \emptyset$. In particular, $P$ is simple if and only if $P^\Delta$ is simplicial, in which case $\partial P^\Delta$ is isomorphic, as a simplicial complex, to the naive dual complex $\{J\subseteq I : \bigcap_{j\in J}D_j\neq \emptyset\}$. 

To define the evaluation map in the context of polytopes, we choose $\delta>0$ small, and define the sets $N_i(P)$ and $N^*(P)$, and the partition of unity $\{\varphi_i^P\}_{i\in I}$ as in Section \ref{sec:evaluation map}. We can then define 
\[
    \ev^P\colon N^*(P)\to \Delta(P)=\partial P^\Delta,\;\quad x\mapsto \sum_{i\in I}\varphi_i^P(x) v^\Delta_i \ .
\]
As for all $J\subseteq I$, we have that  $\bigcap_{j\in J} N_i(P)\neq \emptyset$ if and only if $\{v_j^{\Delta}: j\in J\}$ lie in a common facet, the map $\ev^P$ is well-defined.

Applying the same construction to $P^\Delta$, that is, choosing functions $\beta^{\Delta}_k=| l^{\Delta}_k|$ for every facet $D_k^\Delta$, $k\in K$, of $P^\Delta$, we obtain an evaluation map
\[
    \ev^{P^\Delta}\colon N^*(P^\Delta)\to \Delta(P^\Delta)=\partial P, \;\quad x\mapsto \sum_{k \in K}\varphi_k^{P^\Delta}(x)v_k\ ,
\]
where $v_k$ is the vertex of $P$ corresponding to the facet $D_k^\Delta$ for $k\in K$.

The connection to the algebraic setup in Section \ref{sec:evaluation map} is through toric varieties. A lattice polytope $P$ defines a projective toric variety $X=X_P$, which is a compactification of a big open torus $S$ acting on $X$. The $S$-action on $X$ yields an algebraic moment map
\[
    \mu\colon X\to P
\]
with $\mu^{-1}\mathring P= S$ the big torus. Moreover, the orbit closures in $X$ are precisely the preimages of the faces of $P$. 

If $X$ is smooth, the polytope $P$ is simple. 
As $\mu$ is semialgebraic, we may choose rug functions $\beta^P_i$ for the polytope and define the rug functions for $X$ as $\beta^X_i=\beta^P_i\circ \mu$. In particular, after identifying $\partial P^\Delta$ with $\Delta(X)$ in the canonical way, we have
\[
\ev^X \equ \ev^P\circ \mu\ .
\]
\begin{proposition}
\label{prop:linear homotopies work}
Let $P$ be as above (not necessarily simple). For $0\ll \epsilon<1$, the map
\[
    \widetilde \ev^P\colon N^*(P)\to N^*(P^\Delta)\ ,\;\quad  x\mapsto \epsilon \ev^{P}(x)
\]
is well-defined. Moreover, if $\epsilon$ is sufficiently close to $1$, then for every $x\in N^*(P)$, the straight line between $x$ and $\widetilde\ev^{P^\Delta}(\widetilde\ev^P(x))$
is contained in $N^*(P)$.
\end{proposition}
\begin{proof}
The fact that $\widetilde \ev^P$ is well-defined for $\epsilon$ close to $1$  follows directly from the compactness of $\partial P^\Delta$.    
Now fix $\epsilon<1$ sufficiently close to $1$, and let $x\in N^*(P)$.
By our choice of rug functions, a point $y$ in the interior of $P$ lies in $N_i(P)$ if and only if $\langle v^\Delta_i,y\rangle>1-\delta$, and a point $m$ in the interior of $P^\Delta$ lies in $N_k(P^\Delta)$ if and only if $\langle m,v_k\rangle>1-\delta$.
Since $\langle v_k,v^\Delta_i\rangle\leq1$, with equality if and only if $v_k\in D_i$, there is a constant $0\leq c<1$ such that $\langle v_k,v^\Delta_i\rangle\leq c$ whenever $v_k\notin D_i$. This condition is symmetric in $P$ and $P^\Delta$, since $v_k\in D_i$ if and only if $v^\Delta_i\in D^\Delta_k$, and we assume that $\delta<(1-c)/\max\{|I|,|K|\}$. 
Set 
\begin{align*}
    S&=\{i\in I : x\in N_i(P)\} \quad \text{ and}\\
    T &=\{k\in K:v_k\in \bigcap_{i\in S} D_i\}\ .
\end{align*}
Then $\ev^P(x)$ is in the relative interior of $F= \conv\{v_i^\Delta : i\in S\}$, which is contained in $\bigcap_{k\in T}D^\Delta_k$. 
Because $\epsilon$ was chosen sufficiently close to $1$,  we have 
\[
\widetilde \ev^P(x)=\epsilon \ev^P(x)\in \bigcap_{k\in T} N_k(P^\Delta)\ . 
\]
Therefore, $\ev^{P^\Delta}(\widetilde\ev^P(x))$, which is a convex combination of the vertices $v_k$ with $\widetilde\ev^P(x)\in N_k(P^\Delta)$, is contained in the smallest face of $P$ containing these vertices. This face contains
\[
    \conv\{v_k:k\in T\}=\bigcap_{i\in S}D_i\ .
\]
Let $S'$ be the set of all $i\in S$ such that $D_i$ contains this face, so that $\ev^{P^\Delta}(\widetilde\ev^P(x))\in\bigcap_{i\in S'}D_i$.
The set $S'$ is non-empty: in fact, since the $\varphi^P_i(x)$ sum up to $1$ and vanish for $i\notin S$, there is some $i_0\in S$ with $\varphi^P_{i_0}(x)\geq1/|I|$. If $v_k\notin D_{i_0}$, then
\[
    \big\langle v_k,\widetilde\ev^P(x)\big\rangle
    =\epsilon\sum_{i\in S}\varphi^P_i(x)\,\langle v_k,v^\Delta_i\rangle
    \leq1-\varphi^P_{i_0}(x)(1-c)
    \leq1-\frac{1-c}{|I|}<1-\delta\ ,
\]
so that $\widetilde\ev^P(x)\notin N_k(P^\Delta)$. Hence every vertex $v_k$ with $\widetilde\ev^P(x)\in N_k(P^\Delta)$ lies in $D_{i_0}$, and $i_0\in S'$.
Again by the choice of $\epsilon$, the point
\[
    \widetilde\ev^{P^\Delta}(\widetilde \ev^P(x))=\epsilon \ev^{P^\Delta}(\widetilde\ev^P(x))
\]
is contained in $\bigcap_{i\in S'}N_i(P)$. As $x$ is also contained in $\bigcap_{i\in S'}N_i(P)$ and the set $\bigcap_{i\in S'}N_i(P)\setminus\partial P$ is convex and contained in $N^*(P)$, it follows that the line segment between $x$ and $\widetilde \ev^{P^\Delta}(\widetilde\ev^P(x))$ is fully contained in $N^*(P)$.
\end{proof}

\begin{corollary}
\label{cor:evaluation map for polytopes homotopy}
The evaluation map 
\[
    \ev^P\colon N^*(P)\to \partial P^\Delta
\]
is a homotopy equivalence. Moreover, if $H$ is a finite group acting linearly on $K_\R$ and preserving $P$, then $\ev^P$ is an $H$-equivariant  homotopy equivalence. 
\end{corollary}

\begin{proof}
Let 
\begin{align*}
\widetilde \ev^P\colon&
     N^*(P)\to N^*(P^\Delta) \quad\text{ and}\\
\widetilde \ev^{P^\Delta}\colon&
     N^*(P^\Delta)\to N^*(P)
\end{align*}
be the maps from Proposition \ref{prop:linear homotopies work} for $\epsilon$ close enough to $1$. Then by Proposition \ref{prop:linear homotopies work}, the composites $\widetilde\ev^P\circ \widetilde\ev^{P^\Delta}$ and $\widetilde\ev^{P^\Delta}\circ\widetilde \ev^P$ are connected to the respective identities via straight-line homotopies. If a finite group acts linearly and both $\ev^P$ and $\ev^{ P^\Delta}$ are  $H$-equivariant, then the linear homotopy is $H$-equivariant as well. So in this case, $\widetilde \ev^P$ is an $H$-equivariant homotopy equivalence.

By definition, $\widetilde \ev^P$ is the composite
\[
    N^*(P)\xrightarrow{\ev^P}\partial P^\Delta
        \xrightarrow{\cdot\epsilon}N^*(P^\Delta)\ .
\]
The second map is a homotopy equivalence, because one can scale the points of $N^*(P^\Delta)$ to a unique point in $\epsilon\partial P^\Delta$, and the scaling is $H$-equivariant if $H$ acts linearly on $P$. It follows that $\ev^P$ is a homotopy equivalence as well, $H$-equivariantly if $\ev^P$ is $H$-equivariant. 
\end{proof}

The above considerations lead to  the following toric version of the geometric $P=W$ conjecture.

\begin{theorem}
\label{thm:toric p=w}
Let $X$ be a smooth projective toric variety compactifying a big open torus $S$ with cocharacter lattice $K$. Then there is a homotopy equivalence 
\(
   \Delta(\partial X) \xrightarrow{\chi}K_\R\setminus\{0\}
\)
making the triangle
\[
\begin{tikzcd}
    &
        N^*(S)
        \arrow[dl,swap, "\ev^S"]
        \arrow[dr, "-\log|\, \cdot\, |"] &
            \\
   \Delta(\partial X) \arrow[rr,"\chi"]&
        &
            K_\R\setminus \{0\} \ .
\end{tikzcd}
\]
commute up to homotopy. Moreover, if a finite group $H$ acts on $X$ via toric isomorphisms, then $\ev^S$ can be chosen to be $H$-equivariant and the triangle commutes up to an $H$-equivariant homotopy.
\end{theorem}

\begin{proof}
Let $P$ be a lattice polytope with associated toric variety $X=X_P$. After scaling and shifting the polytope $P$ we may assume that $0$ is its center of mass.
If a finite group $H$ acts on $X$ via toric isomorphisms, the normal fan of $P$ is $H$-invariant. Consequently, replacing $P$ with the Minkowski sum $P' \coloneqq \sum_{h\in H} hP$ does not change its normal fan, so that $X_{P'} \cong X_P = X$. We may thus assume without loss of generality that $P$ is invariant under $H$.   Let $\mu\colon X_P\to P$ be the algebraic moment map. Choose affine rug functions $\beta_i^P$ for the facets of $P$ as in Section \ref{sec:evaluation map polytopes}. Summing over all shifts as in \eqref{eq:sum of shifts} makes sure that $N^*(P)$ is $H$-invariant and the evaluation map $\ev^P\colon N^*(P)\to \Delta(\partial X)$ is $H$-equivariant. Define the rug functions for $X$ by pulling back along $\mu$, that is $\beta^X_i=\beta^P_i\circ \mu$. As $\mu$ is $H$-equivariant, each $N_i(X)$ is $H$-invariant and the evaluation map $\ev^X=\ev^P\circ \mu$ is $H$-equivariant as well. 

The moment map $\mu(x)$ only depends on the absolute values of $x\in S$. The restriction $\mu|_S$  factors through $-\log|\, \cdot\, |$, that is, there exists an $H$-equivariant map $\widetilde \mu\colon K_\R\to \mathring P$ with $\mu|_S=\widetilde \mu\circ (-\log|\, \cdot \, |)$. In fact, $\widetilde \mu$ is a homeomorphism \cite[Theorem 12.2.5]{CLS11}. Thus, we  have the  diagram
\[
\begin{tikzcd}
    &
        N^*(X)
        \arrow[dl,swap,"-\log|\, \cdot\, |"]
        \arrow[d,"\mu"]
        \arrow[dr,"\ev^X"]&
            \\
    \widetilde\mu^{-1}N^*(P)
    \arrow[d,"j"]
    \arrow[r,"\widetilde \mu"]&
        N^*(P)
        \arrow[d,"\simeq"]
        \arrow[r,"\ev^P"]&
            \Delta(\partial X)\\
    K_\R\setminus \{0\}
    \arrow[r,"\widetilde\mu"]&
        \mathring P\setminus \{0\} \ ,
       % &   \\\ ,
\end{tikzcd}
\]
where all maps  are $H$-equivariant. 

By Corollary \ref{cor:evaluation map for polytopes homotopy}, the evaluation map $\ev^P$ is a homotopy equivalence $H$-equivariantly. The other   two horizontal maps, restrictions of $\widetilde \mu$, are homeomorphisms. The inclusion of $N^*(P)$ into $\mathring P\setminus\{0\}$ is a homotopy equivalence $H$-equivariantly, as one can deformation retract both sets linearly to the even smaller set $\omega\partial P$ for some $\omega<1$ close to $1$. Because $\widetilde \mu$ is a homeomorphism, the inclusion $j\colon\widetilde{\mu}^{-1}N^*(P)\rightarrow K_{\mathbb{R}}\setminus \{0\}$ is a homotopy equivalence $H$-equivariantly as well. 
%The assertion now follows by taking for $\chi$ the composite of an $H$-equivariant homotopy inverse of the inclusion $\widetilde \mu^{-1}N^*(P)\hookrightarrow K_\R\{0\}$ composed with $\ev^P\circ \widetilde \mu$. 
The assertion now follows by defining $\chi \colon \Delta(\partial X) \to K_{\mathbb{R}} \setminus \{0\}$ as the composite of $j\circ \tilde{\mu}^{-1}$ with a homotopy inverse of $\mathrm{ev}^P$. By construction, $\chi$ is an $H$-equivariant homotopy equivalence making the diagram commute up to $H$-equivariant homotopy. 
\end{proof}

\subsection{Finite quotients of toric varieties}

\begin{lemma}
\label{lem:dlt quotient}
Let $K$ be a lattice, and let $H$ be a finite group acting on $K$ linearly. Let $\Sigma$ be a smooth complete fan in $K_\R$ on which $H$ induces an action. Suppose that no element in $H$ acts as a pseudoreflection on $K$. Suppose further that:
\begin{enumerate}[(a)]
\item for every $\sigma\in \Sigma$ and $h\in H$ an equality $h\sigma=\sigma$ implies that $h|_\sigma=\id_\sigma$;
\item for any cone $\sigma \in \Sigma$ and face $\tau \preceq \sigma$, if $h\tau \preceq \sigma$ for some $h \in H$, then $h\tau = \tau$.
\end{enumerate}

Let $X=X_\Sigma$ be the toric variety associated to $\Sigma$. Then $(X/H, \partial X/H)$ is a dlt pair.
\end{lemma}

\begin{proof}
Without loss of generality, we may assume that $H$ acts faithfully. Let $\sigma\in \Sigma$ and let $O(\sigma)$ denote the associated torus orbit. Suppose $h\in H$ stabilizes $O(\sigma)$, which is equivalent to $h\sigma=\sigma$. By assumption (a), this implies that $h$ stabilizes $\sigma$ pointwise. Assume $h$ acts trivially on $O(\sigma)$. This is equivalent to $h$ acting trivially on $K_\R/\Span(\sigma)$. Because $h$ also acts trivially on $\Span(\sigma)$, it acts unipotently on $K_\R$. But it also has finite order, forcing it to act trivially. Because the action of $H$ is faithful, we conclude that $h=1$. We conclude that the torus orbits have generically trivial stabilizers.

As every nontrivial $h\in H$ is neither a pseudoreflection, nor fixes a boundary divisor pointwise, the action of $H$ is free in codimension $1$. Therefore, the pair $(X/H,\partial X/H)$ is log canonical by \cite[Proposition 5.20]{KM98}. The log canonical centers of $(X/H,\partial X/H)$ are precisely the images of the log canonical centers of $(X,\partial X)$. The log canonical centers of $(X, \partial X)$, in turn, are precisely the orbit closures \cite[Section 3.9]{Fujino}. So let $Z$ be an orbit closure. We have seen that the points of $Z$ have generically trivial stabilizers, that is, $Z$ intersects the locus $X^{\mathrm{free}}$, where the action of $H$ is free. The quotient map $X^{\mathrm{free}} \to X^{\mathrm{free}}/H$ is étale, making $X^{\mathrm{free}}/H$ smooth; moreover, hypothesis~(b) rules out self-intersections of the boundary divisor $\partial X/H$. Consequently, $X^{\mathrm{free}}/H$ is contained in the snc locus of $(X/H, \partial X/H)$. In particular, the image of $Z$ in $X/H$ intersects the snc locus non-trivially. This shows that no lc center of $(X/H,\partial X/H)$ is contained in the complement of the snc locus, which is precisely the dlt condition. 
\end{proof}

\begin{proposition}
\label{prop:barycentric subdivision}
Let $K$ be a lattice, let $H$ be a finite group acting linearly on $K$, and let $\Sigma$ be a smooth fan in $K_\R$ which comes with  an induced $H$-action. Let $\Delta$ denote the barycentric subdivision of $\Sigma$.
\begin{enumerate}[(a)]
    \item If $\delta\in \Delta$, $\tau\preceq \delta$ is a face of $\delta$, and $h\in H$ such that $h\tau$ is also a face of $\delta$, then $h\tau=\tau$ and $h|_\tau=\id_\tau$. In particular, for all $h\in H$ and $\delta\in \Delta$ the equality $h\delta=\delta$ implies $h|_\delta=\id_\delta$.
    \item Suppose that $\Sigma$ already has the property of $\Delta$ from part (a), that is $\tau\preceq \sigma$ and $h\tau\preceq \sigma$ implies $h|_\tau=\id_\tau$.
    Let $\sigma_1,\ldots , \sigma_n$ be rays of $\Delta$ and $h_1,\ldots, h_n\in H$. If both $\sigma_1,\ldots,\sigma_n$ and $h_1\sigma_1,\ldots,h_n\sigma_n$ span a cone of $\Delta$, then there exists $h\in H$ with $h\sigma_i=h_i\sigma_i$.
\end{enumerate}
\end{proposition}

\begin{proof}
We note that the poset of cones in $\Delta$ is isomorphic to the poset of chains of nontrivial cones of $\Sigma$. In particular, rays of $\Delta$ correspond to nontrivial cones of $\Sigma$.

For part (a), let $\mc F$ and $\mc G$ be chains of nontrivial cones of $\Sigma$ corresponding to $\tau$ and $\delta$, respectively. The conditions $\tau\preceq \delta$ and $h\tau\preceq\delta$ are equivalent to $\mc F\subseteq \mc G$ and $h\mc F\subseteq \mc G$. For $\sigma\in \mc F$, both $\sigma$ and $h\sigma\in h\mc F\subseteq \mc G$ are contained in $\mc G$. Because a chain of cones of $\Sigma$ can contain at most one element of dimension $\dim\sigma$, it follows that $\sigma=h\sigma$. As $\sigma$ was an arbitrary element in $\mc F$, corresponding to a ray of $\tau$, it follows that $h$ is the identity on all rays of $\tau$ (note that $h$ is an isomorphism on lattices and thus needs to map primitive ray generators to primitive ray generators) and hence $h|_\tau=\id_\tau$. Taking $\tau=\delta$, we conclude that $h\delta=\delta$ implies $h|_\delta=\id_\delta$.

In part (b), we may identify the rays $\sigma_i$ with nontrivial cones of $\Sigma$. The condition on spanning cones means that both $\{\sigma_i:1\leq i\leq n\}$ and $\{h_i\sigma_i:1\leq i\leq n\}$ are chains in $\Sigma$. After potentially relabeling, we may assume that 
\[
    \sigma_1\subseteq \sigma_2\subseteq \cdots \subseteq \sigma_n \ .
\]
For dimension reasons, this implies that
\[
    h_1\sigma_1\subseteq h_2\sigma_2\subseteq \cdots\subseteq h_n\sigma_n\ .
\]
Acting on the cones $h_i\sigma_i$ by $h_n^{-1}$, we may assume that $h_n=1$. In this case, we have, for $1\leq i\leq n$, that both $\sigma_i$ and $h_i\sigma_i$ are faces of $\sigma_n$. By our assumption on $\Sigma$, it follows that $h\sigma_i=\sigma_i$ for all $i$, so that we can take $h=h_n$.
\end{proof}

\begin{lemma}
\label{lem:boundary structure of quotient}
Let $\Delta$ be a fan as in Proposition \ref{prop:barycentric subdivision} (b), let $X=X_\Delta$ be the associated toric variety, and let $Y=X/H$ and $\pi\colon X\to Y$ the quotient map.  
\begin{enumerate} [(a)]
    \item For any collection of boundary divisors $D_1,\ldots, D_k$ of  $X$, 
    the intersection $\bigcap_{i=1}^k \pi(D_i)$ is irreducible. Moreover, if $\bigcap_{i=1}^k D_i\neq \emptyset$  then 
    \[
        \bigcap_{i=1}^k \pi(D_i)=\pi\left(\bigcap_{i=1}^k D_i\right)\ .
    \]
    \item For two distinct boundary divisors $D$ and $E$ of $X$ with $\pi(D)=\pi(E)$, we have $D\cap E=\emptyset$.
\end{enumerate}
\end{lemma}

\begin{proof}
Let $D_i=V(\rho_i)$ be the closed stratum corresponding to a ray $\rho_i$ of $\Delta$. We have
\[
    \pi^{-1}\left(\bigcap_{i=1}^k \pi(D_i)\right)
    = \bigcup_{h_1,\ldots, h_n\in H} \bigcap_{i=1}^k V(h_i\rho_i)
    = \bigcup_{h\in H} h\left(\bigcap_{i=1}^k D_i\right)\ ,
\]
If $\bigcap_{i=1}^k\pi(D_i)\neq \emptyset$, we may assume, without loss of generality, that $\bigcap_{i=1}^k D_i\neq \emptyset$. It then follows from Proposition \ref{prop:barycentric subdivision} (b) that for any choice of $h_1,\ldots,h_k\in H$ with $\bigcap_{i=1} k V(h_i\rho_i)\neq \emptyset$, this intersection is equal to $h\left(\bigcap_{i=1}^k D_i\right)$ for some choice of $h\in H$. Since $\bigcap_{i=1}^k D_i$ is irreducible, the assertion follows.

For part (b), write $D=V(\rho)$ and $E=V(h\rho)$ for a ray $\rho$ of $\Delta$ and $h\in H$. Assume $D\cap E\neq \emptyset$. Then $\rho$ and $h\rho$ span a two-dimensional cone $\sigma$. But then Proposition \ref{prop:barycentric subdivision} (a) implies that $\rho=h\rho$, contradicting $E\neq D$.
\end{proof}

\subsection{Proof of Theorem \ref{thm:geometric P=W}}

Let $\Gamma=\pi_1(X)^\ab_\tf$ be the torsion free part of the abelianization of $\pi_1(X)$. The moduli space 
\[
    M^T_\Betti(X)\cong M^{\G_m}_\Betti(X)\otimes_\Z N =\Hom(\Gamma,\C^*)\otimes_\Z N\ 
\]
is an algebraic torus. It follows from the description of the non-abelian Hodge correspondence in rank one (see \cite[p.\ 21]{Simpson_nonabelianHodge} for the curve case; the general case is analogous) that the composite $\sd^T_X\circ \eta^T_X$ is the composite of
\[
    -\log|.|\colon \Hom(\Gamma,\C^*)\to \Hom(\Gamma,\R)
\]
with a $\mathbb{R}-$linear isomorphism $\Hom(\Gamma,\R)\cong H^0(X,\Omega^1_X)$. Here, $B^T_X= H^0(X,\Omega^1_X)\otimes_\mathbb{Z}N$ is the Hitchin base.

Let $Z$ be a projective toric variety with big open torus $S=M^T_\Betti(X)$ and fan $\Sigma$. Replacing $\Sigma$ with the common refinement of the fans $w\Sigma$, $w\in W$, that is, with the fan whose cones are the intersections $\bigcap_{w\in W}w\sigma_w$ with $\sigma_w\in\Sigma$, we may assume that $\Sigma$ is $W$-invariant and $W$ acts on $Z$ via toric automorphisms. Equivariantly resolving singularities, we may also assume that $Z$ is smooth. Finally we subdivide $\Sigma$ barycentrically twice. Let $\sum_{i\in I}D_i$ be the boundary of $Z$. By Theorem \ref{thm:toric p=w}, there is a $W$-invariant choice of sets $N^*_i(Z)$ and rug functions compatible with the $W$-actions such that the evaluation map $\ev^S\colon N^*(S)\to \Delta(\partial Z)$ is $W$-equivariant. Moreover, there is a $W$-equivariant homotopy equivalence $\chi\colon \Delta(\partial Z)\to B^T_X\setminus \{0\}$ making the diagram

\begin{equation}
\label{eq:T-case of geometric P=W}
\begin{tikzcd}
    N^T_\Betti(X)^*
    \arrow[d, "\ev^S"]
    \arrow[r, "\eta^T_X"] &
        N^T_\Dol(X)^*
        \arrow[d,"\sd^T_X"]\\
    \Delta(\partial Z)
    \arrow[r,"\chi"]&
        B^T_X\setminus \{0\}\ .
\end{tikzcd}    
\end{equation}

commute up to $W$-equivariant homotopy, where we denote $N^T_\Betti(X)^*= N^*(S)$ and $N^T_\Dol(X)^*=(\sd_X^T)^{-1} (B^T_X \setminus\{0\})$.

Denote $Y=Z/W$, and let $\pi\colon Z\to Y$ be the quotient map. This is a compactification of $S/W=\widetilde M^G_\Betti(X)$, but typically not snc. However since $\dim S\geq 2$, it is a dlt compactification by Proposition \ref{prop:barycentric subdivision} (a) and Lemma \ref{lem:dlt quotient}. Therefore, the dual complex of $(Y,Y\setminus (S/W))$ still makes sense and it is the boundary complex of $\widetilde M^G_\Betti(X)$. Intersections of boundary divisors are irreducible on $Y$ by Proposition \ref{prop:barycentric subdivision} (b) and Lemma \ref{lem:boundary structure of quotient} (a), so condition \labelcref{eq:irrecubile condition} is satisfied, and an evaluation map can be defined exactly as in the snc case. Moreover, it follows from Proposition \ref{prop:barycentric subdivision} that 
\[
\Delta(\partial Y) \equ \Delta(\partial Z)/W\ .
\]
 
Next, we show  how the rug functions on $Z$ induce rug functions on $Y$. For a boundary divisor $D$ of $Z/W$, let $D_j$, $j\in J\subseteq I$ be the boundary divisors of $Z$ lying over it. By Lemma \ref{lem:boundary structure of quotient} (b), the divisors $D_j$, $j\in J$ are pairwise disjoint, and the same is thus true for the neighborhoods $N_j(Z)$. Define the function
\[
\widetilde \beta^Y_D \colon Z\to \R_{\geq 0},\; x\mapsto 
\begin{cases}
    \beta_j^{Z}(x) &\text{if }x\in N_j(Z) \text{ for some }j\in J\\
    \delta &\text{else}\ ,
\end{cases}
\]
where $\delta$ is the small constant we had chosen to define the neighborhoods $N_j(Z)$, that is, 
\[
N_j(Z) \equ \{x:\beta^Z_j(x)\leq \delta\}\ .
\]
As 
\[
\beta^{Z}_j(x) \equ \beta^{Z}_{wj}(wx)
\]
for all $j\in J$, $x\in Z$, and $w\in W$, the function $\widetilde\beta^Y_D$ is $W$-invariant and thus descends to a function $\beta^Y_D$, which is a rug function for $D$. By definition, we have $N_D(Y)=\pi(N_j(Z))$ for every $j\in J$, and if we define the partition of unity for $Y$ exactly as in \eqref{eq:partition of unity}, the evaluation function $\ev^{S/W}$ is precisely the map that the $W$-equivariant map $\ev^{S}$ descends to.

We can now take the quotient of diagram \eqref{eq:T-case of geometric P=W}.  Let us  denote 
\[
N^G_\Betti(X)^* \equ N^*(\widetilde M^G_\Betti(X))\ ,
\]
and 
\[
N^G_\Dol(X)^* \equ  (\widetilde \sd_X^G)^{-1}(\widetilde B^G_X\setminus\{0\})\ ,
\]
where $\widetilde B^G_X= B^T_X/W$ is the Hitchin base. Then we obtain the  diagram
\[
\begin{tikzcd}
    N^G_\Betti(X)^*
    \arrow[r,"="]
    \arrow[rrr,bend left=20,"\eta^G_X"]
    \arrow[d,"\ev^{S/W}"]&[1ex]
        N^T_\Betti(X)^*/W
        \arrow[d, "\ev^S/W"]
        \arrow[r, "\eta^T_X/W"] &[1ex]
            N^T_\Dol(X)^*/W
            \arrow[d,"\sd^T_X/W"]
            \arrow[r,"="]&[1ex]
                N^G_\Dol(X)^*
                \arrow[d,"\widetilde\sd^G_X"]\\[3ex]
    \Delta(\partial Y)
    \arrow[r,"="]&
        \Delta(\partial Z)/W
        \arrow[r,"\chi/W"]&
            (B^T_X\setminus \{0\})/W
            \arrow[r,"="]&
                \widetilde B^G_X\setminus\{0\}\ ,
\end{tikzcd}    
\]
which commutes up to homotopy. Since  $\chi$ was a $W$-equivariant homotopy equivalence, we are done.

%%%%%%%%%%%%%%%%%%%%%%%%%%%%%%%%%%%%%%%%%%%%%%%%%%%
%% SECTION 5
%%%%%%%%%%%%%%%%%%%%%%%%%%%%%%%%%%%%%%%%%%%%%%%%%%%

\section{Explicit $P=W$ for classical groups and complements}
\label{sec:pwclassicalgroups}

To complement the material in the previous sections, we present a detailed description of the $P=W$ picture (including the moduli spaces in question and the filtrations on the cohomology) for the classical groups of type $C_n$ and $D_n$. The case $A_n$ has been treated in detail in \cite{BKU23}. In particular, when $X$ is an elliptic curve, the moduli spaces of $G$-Higgs bundles for the classical groups and the corresponding Hitchin maps have been described explicitly in \cite{FGN14}.

The case of type $B_n$ is very similar to $C_n$. 
Throughout the whole section, thanks to Proposition \ref{prop:reduction to T Betti} we can assume that $X$ be a smooth projective variety with free abelian group $\pi_1(X)\cong \mathbb{Z}^{2g}$.

\subsection{Type $C_n$: Symplectic Groups $\Sp_{2n}$}
\label{subsec:type-Cn}

Let us specialize to the symplectic group $G = \Sp_{2n}$. 
A maximal torus $T \subseteq \Sp_{2n}$ is isomorphic to $(\mathbb{C}^*)^n$, with cocharacter lattice $N \cong \mathbb{Z}^n$ and Lie algebra $\mathfrak{t} \cong \mathbb{C}^n$. The Weyl group is the hyperoctahedral group
\[
W(C_n) = S_n \ltimes (\mathbb{Z}_2)^n\ ,
\]
acting on $\mathfrak{t}$ and $N$ by permutations and sign changes of coordinates. The invariant ring $\mathbb{C}[\mathfrak{t}]^{W(C_n)}$ is freely generated by the elementary symmetric polynomials in the squared coordinates:
\[
f_k(z_1, \dots, z_n) = e_k(z_1^2, \dots, z_n^2), \quad k = 1, \dots, n\ ,
\]
of respective  degrees $2, 4, \dots, 2n$.

\subsubsection{Dolbeault moduli space and spectral data}
Consider now the varieties
\[
\mathrm{Pic}^0(X) \otimes_\mathbb{Z} N \cong (\mathrm{Pic}^0(X))^n, \quad \text{and}\quad H^0(X, \Omega_X^1)\otimes_\Z N \cong (\mathbb{C}^g)^n\ .
\]
Let $S \coloneqq (\mathrm{Pic}^0(X) \times \mathbb{C}^g)/\mathbb{Z}_2$ be the quotient under the diagonal sign involution $[-1] \times (-1)$. Commuting the quotients identifies the Dolbeault moduli space with the $n$-th symmetric product of $S$:
\[
M_{\Dol}^{\Sp_{2n}}(X)\cong 
\frac{(\mathrm{Pic}^0(X) \times \mathbb{C}^g)^n}{S_n \ltimes (\mathbb{Z}_2)^n} \cong 
\frac{S^n}{S_n}= 
S^{(n)}\ .
\]

Under this identification, the spectral data morphism specializes to 
\[
\sigma \colon S^{(n)} \longrightarrow (\mathbb{C}^g / \mathbb{Z}_2)^{(n)}, \quad \{ [L_i, z_i] \}_{i=1}^n \longmapsto \{ [z_i] \}_{i=1}^n
\]
and the Hitchin map $h$ factors as $h=\iota\circ\sigma$, where
\[
    \iota\colon(\C^g/\Z_2)^{(n)}\longrightarrow\bigoplus_{k=1}^nH^0(X,\Sym^{2k}\Omega^1_X)
\]
is induced by the polarizations of the invariant generators $e_k(t_1^2,\dots,t_n^2)$ of $\C[\mathfrak t]^{W(C_n)}$. More precisely, write a point of $(\C^g)^n=H^0(X,\Omega^1_X)\otimes_\Z N$ as $(z_1,\dots,z_n)$ with $z_i\in H^0(X,\Omega^1_X)$. Polarizing $e_k(t_1^2,\dots,t_n^2)$ amounts to evaluating it at $(z_1,\dots,z_n)$ in the commutative ring $\Sym^\bullet H^0(X,\Omega^1_X)$, followed by the natural map $\Sym^{2k}H^0(X,\Omega^1_X)\to H^0(X,\Sym^{2k}\Omega^1_X)$. Thus $\iota$ sends $\{[z_i]\}_{i=1}^n$ to the image of $\big(e_1(z_1^2,\dots,z_n^2),\dots,e_n(z_1^2,\dots,z_n^2)\big)$, where $z_i^2\in\Sym^2H^0(X,\Omega^1_X)$. Since this expression only depends on $\pm z_i$ and is symmetric in $i$, it is well-defined on $(\C^g/\Z_2)^{(n)}$, and by Lemma \ref{lem:finiteinjectivebetaxg} the morphism $\iota$ is finite and injective on points.
By Proposition \ref{prop:allperverseequal}, the perverse filtration induced by $h$ coincides with that induced by $\sd^G_X$.

\begin{proposition}\label{prop:Cn-Dolbeault-filtration}
The perverse filtration associated with $\sigma$ on $H^*(M_{\Dol}^{\Sp_{2n}}(X), \mathbb{Q})$ is given by
\[
P_p H^d(M_{\Dol}^{\Sp_{2n}}(X), \mathbb{Q}) = 
\begin{cases} 
0 & \text{if } p < d \\ 
H^d(M_{\Dol}^{\Sp_{2n}}(X), \mathbb{Q}) & \text{if } p \ge d\ . 
\end{cases}
\]
\end{proposition}
\begin{proof}
This is an immediate application of K\"unneth formula and Proposition \ref{prop:quotient-filtration} applied to the product $\left(\operatorname{Pic}^0(X)\right)^n\times(\mathbb{C}^g)^n$, with $W = W(C_n)$ acting diagonally.
\end{proof}

\subsubsection{Betti moduli space and cohomological $P=W$}

The Betti moduli space is given by
\[
    M_{\Betti}^{\Sp_{2n}}(X) \cong 
    \frac{(\C^*)^{2g}\otimes_\Z N}{S_n \ltimes (\mathbb{Z}_2)^n} \cong
    \frac{((\mathbb{C}^*)^{2g})^n}{S_n \ltimes (\mathbb{Z}_2)^n} \cong 
    \frac{Y^n}{S_n}\cong
    Y^{(n)}\ ,
\]
where $Y \coloneqq (\mathbb{C}^*)^{2g}/\mathbb{Z}_2$ is the quotient of the algebraic torus by the inversion involution $t \mapsto t^{-1}$.

\begin{proposition}\label{prop:Cn-Betti-weights}
The cohomology of $M_{\Betti}^{\Sp_{2n}}(X)$ satisfies:
\begin{enumerate}
    \item $H^d(M_{\Betti}^{\Sp_{2n}}(X), \mathbb{Q}) = 0$ for all odd degrees $d$.
    \item For any even degree $d = 2m$, $H^{2m}(M_{\Betti}^{\Sp_{2n}}(X), \mathbb{Q})$ is pure of weight $4m$.
\end{enumerate}
In particular, $W_{2k} H^d = H^d$ for $k \ge d$ and vanishes for $k < d$.
\end{proposition}
\begin{proof}
Note that, under inversion on $(\mathbb{C}^*)^{2g}$, odd cohomology is anti-invariant ($(-1)^* = -\mathrm{id}$ on $H^1$). Hence $H^*(Y, \mathbb{Q}) \cong (H^*((\mathbb{C}^*)^{2g}, \mathbb{Q}))^{\mathbb{Z}_2}$ is concentrated in even degrees, with pure weight $4k$ in degree $2k$. Since the operation of symmetric product is compatible with mixed Hodge structures, the cohomology of  $Y^{(n)}$ inherits this property, that is the isomorphism
\[
H^*(M_{\Betti}^{\Sp_{2n}}(X), \mathbb{Q}) \cong \mathrm{Sym}^n H^*(Y, \mathbb{Q})
\]
is an isomorphism of mixed Hodge structures, which preserves both parity vanishing and weight doubling.
\end{proof}

\begin{corollary}[$P=W$ for $\Sp_{2n}$]\label{cor:PW-Cn}
Under the non-abelian Hodge isomorphism, we have:
\[
P_p H^*(M_{\Dol}^{\Sp_{2n}}(X), \mathbb{Q}) = W_{2p} H^*(M_{\Betti}^{\Sp_{2n}}(X), \mathbb{Q}) = W_{2p+1} H^*(M_{\Betti}^{\Sp_{2n}}(X), \mathbb{Q})\ .
\]
\end{corollary}

\begin{remark} The case of groups of type $B_n$ is analogous to that of type $C_n$, so we do not repeat it. 
\end{remark}

\subsection{Type $D_n$: Even Orthogonal Groups $\SO_{2n}$}
\label{subsec:type-Dn}

We now consider the even orthogonal group $G = \SO_{2n}$ for $n \ge 3$.
A maximal torus $T \subset \SO_{2n}$ is isomorphic to $(\mathbb{C}^*)^n$, with $N \cong \mathbb{Z}^n$ and $\mathfrak{t} \cong \mathbb{C}^n$. The Weyl group $W(D_n)$ is the index $2$ normal subgroup of $W(C_n)$ consisting of signed permutations with an even number of sign changes.
% \[
% W(D_n) = S_n \ltimes (\mathbb{Z}_2)^{n-1} = \ker\Big( \prod_{i=1}^n \epsilon_i \colon W(C_n) \longrightarrow \{\pm 1\} \Big)\ .
% \]
The invariant ring $\mathbb{C}[\mathfrak{t}]^{W(D_n)}$ is freely generated by the elementary symmetric polynomials in the squared coordinates $e_1(z_1^2,\ldots,z_n^2), \dots, e_{n-1}(z_1^2,\ldots, z_n^2)$ of degrees $2, 4, \dots, 2n-2$, together with the \emph{Pfaffian}:
\[
\mathrm{Pf}(z_1, \dots, z_n) \coloneqq z_1 z_2 \cdots z_n\ ,
\]
which has degree $n$ and satisfies $\mathrm{Pf}(z)^2 = e_n(z_1^2, \dots, z_n^2)$.

\subsubsection{The Dolbeault moduli space and spectral data}
The Dolbeault moduli space is given by:
\[
M_{\Dol}^{\SO_{2n}}(X) \cong \frac{(\mathrm{Pic}^0(X) \times \mathbb{C}^g)^n}{W(D_n)}\ .
\]
Since $W(D_n) \subset W(C_n)$ has index $2$, the quotient group $\mathbb{Z}_2 \cong W(C_n)/W(D_n)$, generated by a single sign flip $\tau \colon (L_n, z_n) \mapsto (-L_n, -z_n)$, acts on $M_{\Dol}^{\SO_{2n}}(X)$. This induces a canonical degree $2$ finite morphism (a branched double cover):
\[
\pi \colon M_{\Dol}^{\SO_{2n}}(X) \longrightarrow M_{\Dol}^{\Sp_{2n}}(X) \cong S^{(n)}\ ,
\]
branched over the locus of points $\{[L_i,z_i]\}_{i=1}^n\in S^{(n)}$ such that, for some $i$, the pair $(L_i,z_i)$ is a $2$-torsion point of $\Pic^0(X)\times\C^g$.

The spectral data base $\widetilde{B}^{\SO_{2n}}_X = (\mathbb{C}^g)^n / W(D_n)$ is likewise a branched double cover $\bar{\pi} \colon \widetilde{B}^{\SO_{2n}} _X\to (\mathbb{C}^g/\mathbb{Z}_2)^{(n)}$, and the spectral data morphism $\sd^{\SO_{2n}}_X$ fits into the commutative diagram:
\begin{equation}\label{eq:diagram-Dn-cover}
\begin{tikzcd}[row sep=3em, column sep=4.5em]
M_{\Dol}^{\SO_{2n}}(X) \arrow[r, "\sd_X^{\SO_{2n}}"] \arrow[d, "\pi"']
    & \widetilde{B}^{\SO_{2n}}_X \arrow[d, "\bar{\pi}"] \\
S^{(n)} \arrow[r, "\sigma"']
    & (\mathbb{C}^g/\mathbb{Z}_2)^{(n)}\ .
\end{tikzcd}
\end{equation}

\begin{proposition}\label{prop:Dn-Dolbeault-filtration}
The perverse filtration associated with $\sd^{\SO_{2n}}_X$ on $H^*(M_{\Dol}^{\SO_{2n}}(X), \mathbb{Q})$ is given by:
\[
P_p H^d(M_{\Dol}^{\SO_{2n}}(X), \mathbb{Q}) = 
\begin{cases} 
0 & \text{if } p < d \\ 
H^d(M_{\Dol}^{\SO_{2n}}(X), \mathbb{Q}) & \text{if } p \ge d\ . 
\end{cases}
\]
\end{proposition}
\begin{proof}
Since $\sd^{\SO_{2n}}_X$ is the quotient of the diagonal projection $\mathrm{pr}_2 \colon (\Pic^0(X))^n \times (\mathbb{C}^g)^n \to (\mathbb{C}^g)^n$ under the finite group $W(D_n)$, the K\"unneth formula and Proposition \ref{prop:quotient-filtration} imply that the perverse filtration is governed only by the cohomological degree.
\end{proof}

\subsubsection{Betti moduli space and cohomological $P=W$}
On the Betti side, the character variety is the finite double cover of $Y^{(n)}$:
\[
M_{\Betti}^{\SO_{2n}}(X) \cong \frac{(\mathbb{C}^*)^{2gn}}{W(D_n)} \xrightarrow{\quad \pi_B \quad} Y^{(n)}\ .
\]

\begin{proposition}\label{prop:Dn-Betti-weights}
Every rational cohomology class $\alpha \in H^d(M_{\Betti}^{\SO_{2n}}(X) , \mathbb{Q})$ is pure of weight $2d$. Consequently, Deligne's weight filtration satisfies:
\[
W_{2k} H^d(M_{\Betti}^{\SO_{2n}}(X) , \mathbb{Q}) = 
\begin{cases} 
0 & \text{if } k < d \\ 
H^d(M_{\Betti}^{\SO_{2n}}(X), \mathbb{Q}) & \text{if } k \ge d \ . 
\end{cases}
\]
\end{proposition}
\begin{proof}
The rational cohomology of the algebraic torus $(\mathbb{C}^*)^{2gn}$ is pure of weight $2d$ in degree $d$. Since $M_{\Betti}^{\SO_{2n}}(X) $ is a quotient of $(\mathbb{C}^*)^{2gn}$ by the finite group $W(D_n)$, its cohomology identifies with the $W(D_n)$-invariant sub-mixed Hodge structure:
\[
H^d(M_{\Betti}^{\SO_{2n}}(X) , \mathbb{Q}) \cong \big( H^d((\mathbb{C}^*)^{2gn}, \mathbb{Q}) \big)^{W(D_n)}\ .
\]
Being a sub-mixed Hodge structure of a pure Hodge structure of weight $2d$, it is itself pure of weight $2d$.
\end{proof}

\begin{corollary}[$P=W$ for $\SO_{2n}$]\label{cor:PW-Dn}
Under the non-abelian Hodge correspondence, the $P=W$ identity holds:
\[
P_p H^*(M_{\Dol}^{\SO_{2n}}(X) , \mathbb{Q}) = W_{2p} H^*(M_{\Betti}^{\SO_{2n}}(X) , \mathbb{Q}) = W_{2p+1} H^*(M_{\Betti}^{\SO_{2n}}(X) , \mathbb{Q})\ .
\]
\end{corollary}

\begin{remark}[Adjoint groups]
\label{rem:adjoint groups}
The cohomological statements of this section hold line by line for the adjoint groups $\mathrm{SO}_{2n+1}$, $\mathrm{PSp}_{2n}$ and $\mathrm{PSO}_{2n}$ of respectively type $B_n$, $C_n$ and $D_n$. Indeed, passing to the adjoint group does not change the Weyl group $W$ nor its action on $\mathfrak t$, so the spaces of spectral data, the morphisms $\sigma$ and $\iota$, and the Hitchin map are the same as above. The only change is that the cocharacter lattice $N$ is replaced by the coweight lattice $N_{\mathrm{ad}}\supseteq N$. It is straightforward to see that all the arguments above apply equally to $N_{\mathrm{ad}}$. Hence the perverse filtration is determined only by the cohomological degree, $H^d(M^G_\Betti(X),\Q)$ is pure of weight $2d$, and the cohomological $P=W$ conjecture holds, in accordance with Theorem \ref{thm:P=W}. Moreover, the longest element of $W$ acts as $-1$ on $\mathfrak t$ for types $B_n$, $C_n$, and $D_n$ with $n$ even, so in these cases the odd cohomology vanishes as in Proposition \ref{prop:Cn-Betti-weights}.

What changes is the explicit description of the moduli spaces. For $A=\Pic^0(X)\times\C^g$ or $A=(\C^*)^{2g}$, the natural map $A\otimes_\Z N\to A\otimes_\Z N_{\mathrm{ad}}$ is surjective with finite kernel $\operatorname{Tor}_1^\Z(A,N_{\mathrm{ad}}/N)$, which is given by torsion points of $A$. Hence the moduli spaces for the adjoint groups are quotients of the symmetric products and double covers described above by finite groups of torsion translations.
Also, for $\mathrm{PSp}_{2n}$ and $\mathrm{PSO}_{2n}$ these descriptions hold for the normalizations $\widetilde M^G_\Dol(X)$ and $\widetilde M^G_\Betti(X)$. This does not affect the statements above, since the normalization maps are homeomorphisms and, on the Betti side, induce isomorphisms of mixed Hodge structures.
\end{remark}
%%%%%%%%%%%%%%%%%%%%%%%%%%%%%%%%%%%%%%%%%%%%

\subsection{Symplectic Singularities and Resolutions}
\label{sec:symplectic-resolutions}

In this section, we apply our results from Section 2 to study the singularities of the moduli spaces $M^G_\Dol(X)$ and $M^G_\Betti(X)$, and we show that they admit no symplectic resolution for $g\geq 2$, and we decide the question for $g=1$ in many cases. Throughout the section, $X$ denotes a smooth complex projective variety with nilpotent fundamental group, and $2g=b_1(X)$ denotes the rank of the abelianization of $\pi_1(X)$.

By Proposition \ref{prop:reduction to T Betti} and \Cref{lem:functoriality of spectral data}, the normalizations of these moduli spaces are finite quotients of smooth varieties:
\begin{equation}\label{eq:identification}
\widetilde M^G_\Dol(X)\cong M^T_\Dol(X)/W\cong\big((\Pic^0(X)\times H^0(X,\Omega^1_X))\otimes_\Z N\big)/W\ ,
\end{equation}
and
\[
\widetilde M^G_\Betti(X)\cong M^T_\Betti(X)/W\cong\big((\C^*)^{2g}\otimes_\Z N\big)/W\ .
\]
The above quotients coincide with the corresponding quotients for the Albanese variety of $X$.  
If $G$ is a classical group, then $M^G_\Betti(X)$ is normal: by Proposition \ref{prop:nilpotent-mod-torsion for character varieties}, it only depends on the abelianization of $\pi_1(X)$ modulo torsion, so that we may apply \cite[Theorem 4]{BFN24} to the Albanese variety of $X$. By Simpson's isosingularity principle \cite[Theorem 10.6]{SimII}, the moduli space $M^G_\Dol(X)$ is then normal as well. Since symplectic singularities and symplectic resolutions only make sense for normal varieties, we work with the normalizations throughout; by the above, this is no restriction when $G$ is a classical group.

Recall that a normal variety $M$ has \emph{symplectic singularities} in the sense of Beauville \cite{Be00} if   there is a holomorphic symplectic $2$-form on its smooth locus which can be extended to a holomorphic $2$-form on some --- and thus on every --- resolution of singularities $\pi\colon\widetilde M\to M$. A resolution $\pi$ is called \emph{symplectic} if the extended form is non-degenerate on $\widetilde M$. In this case the resolution is crepant: $K_{\widetilde M}=\pi^*K_M$. Conversely, every crepant resolution of a variety with only symplectic singularities is symplectic.

In view of \eqref{eq:identification}, in order to study the singularities of $M^G_\Dol(X)$ and $M^G_\Betti(X)$, we want to reduce to the case of the maximal torus $T$. To do so, we first define  symplectic forms on $M^T_\Dol(X)$ and $M^T_\Betti(X)$, and then  study the action of the Weyl group on them. In particular, we will show that $M^T_\Dol(X)$ and $M^T_\Betti(X)$ admit translation-invariant holomorphic symplectic forms which are invariant under the action of $W$.

First, we note that the varieties $M^T_\Dol(X)$ and $M^T_\Betti(X)$ are connected commutative algebraic groups, and $W$ acts on them by group automorphisms through its action on $N$.  Fix an ample class $\alpha\in H^2(X,\Q)$ and a $W$-invariant non-degenerate symmetric bilinear form on $\mathfrak t$ (for instance the one induced by a $W$-invariant inner product on $N\otimes_\Z\R$). By the Hard Lefschetz theorem, the alternating form
\[
(a,b)\longmapsto\int_Xa\cup b\cup \alpha^{\dim X-1}
\]
on $H^1(X,\C)$ is non-degenerate, and it induces a perfect pairing between $H^1(X,\mc O_X)$ and $H^0(X,\Omega^1_X)$. The tangent spaces at the identity of $M^T_\Betti(X)$ and $M^T_\Dol(X)$ are
\[
H^1(X,\C)\otimes_\C\mathfrak t
\qquad\text{and}\qquad
\big(H^1(X,\mc O_X)\oplus H^0(X,\Omega^1_X)\big)\otimes_\C\mathfrak t\ ,
\]
respectively, and   the tensor product of the two forms above is a non-degenerate alternating form on both of them. Translating it, we obtain the desired symplectic forms. They are $W$-invariant, because $W$ acts by group automorphisms and preserves the chosen form on $\mathfrak t$.

\begin{theorem}
	\label{thm:symplectic-singularities-general}
	Let $G$ be a connected reductive group. Then $\widetilde M^G_\Dol(X)$ and $\widetilde M^G_\Betti(X)$ are irreducible normal varieties with rational Gorenstein symplectic singularities. In particular, if $G$ is a classical group, the same holds for $M^G_\Dol(X)$ and $M^G_\Betti(X)$.
\end{theorem}

\begin{proof}
	Since $\widetilde M^G_\Dol(X)$ is the quotient of the smooth connected variety $M^T_\Dol(X)$ by $W$, and the group action preserves the holomorphic symplectic form constructed above, we have that $\widetilde M^G_\Dol(X)$ is irreducible and normal, and it has symplectic singularities by \cite[Proposition 2.4]{Be00}. Also, symplectic singularities are rational Gorenstein by \cite[Proposition 1.3]{Be00}. The same argument applies to $\widetilde M^G_\Betti(X)$. The last assertion follows, since for classical groups the moduli spaces in question are normal.
\end{proof}

We now turn to symplectic resolutions. If $G$ is a torus, our  moduli spaces are smooth, so we may assume that $W\neq\{1\}$. The following result shows that the answer depends on $g$ and, for $g=1$, even on the root system of $G$. For abelian character varieties, a very similar result has been proved while this paper was being written, see \cite[Corollary~4.4 and Theorem~4.5]{FGPZ26}. 

\begin{theorem}
	\label{thm:classification-resolutions}
	Let $G$ be a connected reductive group which is not a torus.
	\begin{enumerate}[(i)]
		\item If $g\geq2$, neither $\widetilde M^G_\Dol(X)$ nor $\widetilde M^G_\Betti(X)$ admits a symplectic resolution.
		\item If $g=1$ the following holds. Let $S=(\Pic^0(X)\times\C)/\Z_2$ be the quotient by the involution $-1$, and let $\widetilde S\to S$ be the minimal resolution of its four singular points, which are of type $A_1$. Then $\widetilde M^G_\Dol(X)$ admits a symplectic resolution in the following cases:
		\begin{itemize}
			\item for $G=\GL_n$, the Hilbert--Chow morphism
			\[
			\Hilb^n(\Pic^0(X)\times\C)\longrightarrow(\Pic^0(X)\times\C)^{(n)}\cong\widetilde M^G_\Dol(X)\ ;
			\]
			\item for $G=\SL_n$, the restriction of this morphism over the fiber over $0$ of the summation map $(\Pic^0(X)\times\C)^{(n)}\to\Pic^0(X)\times\C$, which is isomorphic to $\widetilde M^G_\Dol(X)$;
			\item for $G=\Sp_{2n}$ or $G=\SO_{2n+1}$, the composite
			\[
			\Hilb^n(\widetilde S)\longrightarrow\widetilde S^{(n)}\longrightarrow S^{(n)}\cong\widetilde M^G_\Dol(X)\ .
			\]The same holds for $\widetilde M^G_\Betti(X)$, with $\Pic^0(X)\times\C$ replaced by $\C^*\times\C^*$.
		\item If the root system of $G$ has an irreducible component which is not of type $A$, $B$ or $C$ (for instance  $G=\SO_{2n}$ with $n\geq4$, or $G$ is of exceptional type), then neither $\widetilde M^G_\Dol(X)$ nor $\widetilde M^G_\Betti(X)$ admits a projective symplectic resolution.
        \end{itemize}
	\end{enumerate}
\end{theorem}

\begin{proof}
	We treat the Dolbeault case in detail, the Betti case is identical.  For (i), let $w\in W\setminus\{1\}$. Since $w$ acts on $M^T_\Dol(X)$ by a group automorphism, its fixed locus is a closed subgroup with Lie algebra $\big(H^1(X,\mc O_X)\oplus H^0(X,\Omega^1_X)\big)\otimes_\C\mathfrak t^w$. Hence, each of its connected components has codimension $2g\cdot\codim(\mathfrak t^w)\geq2g\geq4$ in $M^T_\Dol(X)$, and the singular locus of $\widetilde M^G_\Dol(X)\cong M^T_\Dol(X)/W$ has codimension at least $4$. 
    
    Observe that the  singular locus  of $\widetilde{M}^G_\Dol(X)$ is non-empty. Indeed, the stabilizer of the identity of $M^T_\Dol(X)$, that is, of the trivial $T$-Higgs bundle, is $W\neq\{1\}$. By the previous computation, no element of $W$ acts on the tangent space at the identity as a pseudoreflection, so that the image of the identity is a singular point by Cartan's linearization theorem and the Chevalley--Shephard--Todd theorem. By Theorem \ref{thm:symplectic-singularities-general} and \cite[Corollary 1]{N01}, the variety $\widetilde M^G_\Dol(X)$ has terminal singularities. As a finite quotient of a smooth variety, it is also $\Q$-factorial \cite[Lemma 5.16]{KM98}. Now let $\pi\colon\widetilde M\to\widetilde M^G_\Dol(X)$ be a symplectic resolution. Since $\widetilde M^G_\Dol(X)$ is $\Q$-factorial, the exceptional locus of $\pi$ is of pure codimension one \cite[Corollary 2.63]{KM98}. Since $\pi$ is crepant and $\widetilde M^G_\Dol(X)$ is terminal, $\pi$ has no exceptional divisors. Hence $\pi$ is finite and birational, and thus an isomorphism by Zariski's Main Theorem. This contradicts the fact that $\widetilde M^G_\Dol(X)$ is singular.
	
	For (ii), note that $M^T_\Dol(X)=(\Pic^0(X)\times\C)\otimes_\Z N$, and that its symplectic form is induced by a translation-invariant symplectic form on $\Pic^0(X)\times\C$. For $G=\GL_n(\C)$, we have $N=\Z^n$, and $W=S_n$ acts by permutations, so that $\widetilde M^G_\Dol(X)\cong(\Pic^0(X)\times\C)^{(n)}$. 
    
    The Hilbert scheme of points on a smooth surface is smooth \cite{F68}, and for a symplectic surface the Hilbert--Chow morphism is a symplectic resolution \cite{B83}. For $G=\SL_n(\C)$, the lattice $N$ is the kernel of the summation map $\Z^n\to\Z$, which splits. Hence $M^T_\Dol(X)$ is the kernel of the summation map $(\Pic^0(X)\times\C)^n\to\Pic^0(X)\times\C$, and $\widetilde M^G_\Dol(X)$ is the fiber over $0$ of the induced map on symmetric products. 
    
    The group $\Pic^0(X)\times\C$ acts on $\Hilb^n(\Pic^0(X)\times\C)$ by translations, and the composite of the Hilbert--Chow morphism with the summation map is equivariant with respect to the action $x\mapsto x+ns$ of an element $s$ on $\Pic^0(X)\times\C$. The latter action is transitive, since $\Pic^0(X)\times\C$ is a divisible group, so all fibers of this composite are isomorphic. As in the construction of generalized Kummer varieties \cite[\S7]{B83}, it follows that its fiber over $0$ is smooth and symplectic, and that it is a symplectic resolution of $\widetilde M^G_\Dol(X)$. 
    
    For $G=\Sp_{2n}$ and $G=\SO_{2n+1}$, we have $N=\Z^n$, and $W=S_n\ltimes (\Z_2)^n$ acts by signed permutations, so that $\widetilde M^G_\Dol(X)\cong S^{(n)}$. The involution $-1$ of $\Pic^0(X)\times\C$ preserves the symplectic form, and its fixed points are the four points of order two. Hence $S$ has four singular points of type $A_1$, and its minimal resolution $\widetilde S\to S$ is crepant. Therefore, the composite $\Hilb^n(\widetilde S)\to\widetilde S^{(n)}\to S^{(n)}$ is a symplectic resolution by \cite{F68,B83}. In the Betti case, the same arguments apply with $\Pic^0(X)\times\C$ replaced by $\C^*\times\C^*$, whose points of order two are again the four fixed points of the involution $-1$.
	
	For the last assertion, the exponential map from the tangent space at the identity to $M^T_\Dol(X)$ is a $W$-equivariant homomorphism of complex Lie groups which is biholomorphic near $0$. Hence a neighborhood of the image of the identity in $\widetilde M^G_\Dol(X)$ is biholomorphic to a neighborhood of the origin in $(\mathfrak t\oplus\mathfrak t)/W\cong(\mathfrak t\oplus\mathfrak t^*)/W$, where $W$ acts diagonally and we identify $\mathfrak t$ with $\mathfrak t^*$ via a $W$-invariant form. A symplectic resolution of $\widetilde M^G_\Dol(X)$ would restrict to a symplectic resolution of this neighborhood.
	
	By \cite{GK04,B09}, if $(\mathfrak t\oplus\mathfrak t^*)/W$ admits a symplectic resolution, then the Calogero--Moser space associated with $W$ is smooth for generic parameters. The latter  only happens if every irreducible component of $W$ is a Weyl group of type $A$, $B$ or $C$. This finishes the proof.
\end{proof}

\begin{example}
	Part (ii) of Theorem \ref{thm:classification-resolutions} applies to varieties of arbitrary dimension. For example, we can consider  $X=E\times F$, where $E$ is an elliptic curve and $F$ is a simply connected smooth projective variety, or  projective bundles over an elliptic curve. In all these cases we have $g=1$, and the normalized moduli spaces of $\GL_n$-, $\SL_n$-, $\Sp_{2n}$- and $\SO_{2n+1}$-Higgs bundles on $X$ admit symplectic resolutions.
\end{example}

\bibliographystyle{plain}
\bibliography{biblio}{}

\end{document}